\documentclass[letterpaper,11pt,reqno]{amsart}

\makeatletter
\usepackage{amssymb}
\usepackage{latexsym}
\usepackage{amsbsy}
\usepackage{amsfonts}
\usepackage{hyperref}
\usepackage{graphicx}

\def\marginpar#1{\ignorespaces}

\DeclareMathOperator\re{Re}
\DeclareMathOperator\leb{Leb}
\DeclareMathOperator\argmax{argmax}
\DeclareMathOperator\argmin{argmin}
\DeclareMathOperator\sgn{sgn}
\DeclareMathOperator\capacity{Cap}
\newtheorem{theorem}{Theorem}[section]
\newtheorem{lemma}[theorem]{Lemma}
\newtheorem{proposition}[theorem]{Proposition}
\newtheorem{corollary}[theorem]{Corollary}

\newtheorem{question}[theorem]{Question}
\newtheorem{openpb}[theorem]{Open problem}
\newtheorem{definition}[theorem]{Definition}
\newtheorem{remark}[theorem]{Remark}
\numberwithin{equation}{section}
\makeatother
\begin{document}
\title[Patterns in random walks and Brownian motion]{Patterns in random walks and Brownian motion}

\author[Jim Pitman]{{Jim} Pitman}
\address{Statistics department, University of California, Berkeley. Email: 
} \email{pitman@stat.berkeley.edu}

\author[Wenpin Tang]{{Wenpin} Tang}
\address{Statistics department, University of California, Berkeley. Email: 
} \email{wenpintang@stat.berkeley.edu}

\date{\today} 

\begin{abstract}
We ask if it is possible to find some particular continuous paths of unit length in linear Brownian motion. Beginning with a discrete version of the problem, we derive the asymptotics of the expected waiting time for several interesting patterns. These suggest corresponding results on the existence/non-existence of continuous paths embedded in Brownian motion. With further effort we are able to prove some of these existence and non-existence results by various stochastic analysis arguments. A list of open problems is presented.
\end{abstract}
\maketitle

\setcounter{tocdepth}{1}
\tableofcontents
\section{Introduction and main results}
\label{intro}
\quad We are interested in the question of embedding some continuous-time stochastic processes $(Z_u, 0 \le u \le1)$ 
into a Brownian path $(B_t; t \geq 0)$, without time-change or scaling, just by a random translation of origin in spacetime. More precisely, we ask the following: 
\begin{question} 
\label{Q1}
Given some distribution of a process $Z$ with continuous paths,
does there exist a random time $T$ such that $(B_{T+u}-B_T; 0 \leq u \leq 1)$ has the same distribution as $(Z_u, 0 \le u \le 1)$?.
\end{question}
\quad The question of whether external randomization is allowed or not, is of no importance here. In fact, we can simply ignore Brownian motion on $[0,1]$, and consider only random times $T \geq 1$. Then $(B_t; 0 \leq t \leq 1)$ provides an independent random element which is adequate for any randomization. 

\quad Note that a continuous-time process whose sample paths have different regularity, e.g. fractional Brownian motion with Hurst parameter $H \ne \frac{1}{2}$, cannot be embedded into Brownian motion. Given $(B_t; t \geq 0)$ standard Brownian motion, we define $g_1:=\sup\{t<1; B_t=0\}$ the time of last exit from $0$ before $t=1$, and $d_1:=\inf\{t>1; B_t=0\}$ the first hitting time of $0$ after $t=1$. The following processes, derived from Brownian motion, are of special interest.
\begin{itemize}
\item
Brownian bridge, which can be defined by
$$\left(b^0_u:=\frac{1}{\sqrt{g_1}}B_{ug_1};0 \leq u \leq 1\right),$$ and its reflected counterpart $(|b^0_u|; 0 \leq u \leq 1)$;
\item
Normalized Brownian excursion defined by
$$\left(e_u:=\frac{1}{\sqrt{d_1-g_1}}|B_{g_1+u(d_1-g_1)}|;0 \leq u \leq 1\right);$$
\item
Brownian meander defined as 
$$\left(m_u:=\frac{1}{\sqrt{1-g_1}}|B_{g_1+u(1-g_1)}|;0 \leq u \leq1\right);$$ 
\item
Brownian co-meander defined as
$$\left(\widetilde{m}_u:=\frac{1}{\sqrt{d_1-1}}|B_{1+u(d_1-1)}|; 0 \leq u \leq 1\right);$$
\item
The three-dimensional Bessel process 
$$\left(R_u:=\sqrt{(B_u)^2+(B'_u)^2+(B^{''}_u)^2}; 0 \leq u \leq 1\right),$$
where $(B'_t; t \geq 0)$ and $(B^{''}_u; u \geq 0)$ are two independent copies of $(B_t; t \geq 0)$;
\item 
The first passage bridge through level $\lambda<0$, defined by
$$(F^{\lambda,br}_u;0 \leq u \leq 1) \stackrel{(d)}{=} (B_u; 0 \leq u \leq 1) \quad \mbox{conditioned on} \quad \tau_{\lambda}=1,$$
where $\tau_{\lambda}:=\inf\{t \geq 0; B_t<\lambda\}$ is the first time at which Brownian motion hits $\lambda<0$.
\item
The Vervaat transform of Brownian motion, defined as 
$$\left(V_u:=  \left\{ \begin{array}{ccl} B_{\tau+u}-B_{\tau}  &\mbox{for} & 0 \leq u \leq 1-\tau \\ B_{\tau-1+u}+B_1-B_{\tau} & \mbox{for} & 1-\tau \leq u \leq 1 \end{array}\right.; 0 \leq u \leq 1\right),$$
where $\tau:=\argmin_{0 \leq t \leq 1}B_t$, and the Vervaat transform of Brownian bridge with negative endpoint $\lambda<0$
$$\left(V_u^{\lambda}:=  \left\{ \begin{array}{ccl} b^{\lambda}_{\tau+u}-b^{\lambda}_{\tau}  &\mbox{for} & 0 \leq u \leq 1-\tau \\ b^{\lambda}_{\tau-1+u}+\lambda-b^{\lambda}_{\tau} & \mbox{for} & 1-\tau \leq u \leq 1 \end{array}\right.; 0 \leq u \leq 1\right),$$
where $(b^{\lambda}_u; 0 \leq u \leq1)$ is Brownian bridge ending at $\lambda \in \mathbb{R}$ and $\tau:=\argmin_{0 \leq t \leq 1}b^{\lambda}_t$.
\end{itemize}
\quad The Brownian bridge, meander, excursion and the three-dimensional Bessel process are well-known. The definition of the co-meander is found in Yen and Yor \cite[Chapter $7$]{YY}. The first passage bridge is studied by Bertoin et al \cite{BCP}. The Vervaat transform of Brownian bridges and of Brownian motion are extensively discussed  in Lupu et al \cite{LPT}. According to the above definitions, the distributions of the Brownian bridge, excursion and (co-)meander can all be achieved in Brownian motion provided some Brownian scaling operation is allowed. Note that the distributions of all these processes are singular with respect to Wiener measure. So it is a non-trivial question whether copies of them can be found in Brownian motion just by a shift of origin in spacetime. Otherwise, for a process $(Z_t, 0 \le t \le 1)$ whose distribution is absolutely continuous with respect to that of $(B_t, 0 \le t \le 1)$, for instance the Brownian motion with drift $Z_t:= \vartheta t + B_t$ for a fixed $\vartheta$, the distribution of $Z$ can easily be obtained as that of $(B_{T+t} - B_T, 0 \le t \le 1)$ for a suitable stopping time $T+1$ by the acceptance-rejection method. We refer readers to Section \ref{33} for further development.

\quad The question raised here has some affinity to the question of embedding a given one-dimensional distribution as the distribution of $B_T$ for a random time $T$.  This {\em Skorokhod embedding problem} traces back  to Skorokhod \cite{Skorokhod} and Dubins \cite{Dubins} -- who found integrable stopping times $T$ such that the distribution of $B_T$ coincides with any prescribed one with mean $0$ and finite second moment. Monroe \cite{Monroe, Monroebis} considered embedding of a continuous-time process into Brownian motion, and showed that every semi-martingale is a time-changed Brownian motion. We refer readers to the excellent survey of Obloj \cite{Obloj}, with a rich body of references. Let $X_t:=(B_{t+u}-B_t;0 \leq u \leq 1)$ for $t \geq 0$ be the moving-window process associated to Brownian motion. In Question \ref{Q1}, we are concerned with the possibility of embedding a given distribution on $\mathcal{C}[0,1]$ as that of $X_T$ for some random time $T$. Indeed, Question \ref{Q1} is the Skorokhod embedding problem for the $\mathcal{C}[0,1]$-valued Markov process $X$.

\quad Let us present the main results of the paper. We start with a list of continuous-time processes that cannot be embedded into Brownian motion by a shift of origin in spacetime.
\begin{theorem}(Non-existence of normalized excursion, reflected bridge, Vervaat transform of Brownian motion, first passage bridge and Vervaat bridge) 
\label{coro1}
\begin{enumerate}
\item
There is no random time $T$ such that $(B_{T+u}-B_T; 0 \leq u \leq 1)$ has the same distribution as $(e_u; 0\leq u \leq 1)$;
\item
There is no random time $T$ such that $(B_{T+u}-B_T; 0 \leq u \leq 1)$ has the same distribution as $(|b_u^0|; 0\leq u \leq 1)$;
\item
There is no random time $T$ such that $(B_{T+u}-B_T; 0 \leq u \leq 1)$ has the same distribution as $(V_u; 0\leq u \leq 1)$;
\item
For each fixed $\lambda<0$, there is no random time $T$ such that $(B_{T+u}-B_T; 0 \leq u \leq 1)$ has the same distribution as $(F^{\lambda,br}_u; 0\leq u \leq 1)$;
\item
For each fixed $\lambda<0$, there is no random time $T$ such that $(B_{T+u}-B_T; 0 \leq u \leq 1)$ has the same distribution as $(V^{\lambda}_u; 0 \leq u \leq 1)$.
\end{enumerate}
\end{theorem}
\quad As we will see later in Theorem \ref{thm2}, Theorem \ref{coro1} is an immediate consequence of the fact that typical paths of normalized excursion, reflected bridge, Vervaat transform of Brownian motion, first passage bridge and Vervaat bridge cannot be achieved in Brownian motion. The next theorem shows the possibility of embedding into Brownian motion some continuous-time processes whose distributions are singular with respect to Wiener measure.
\begin{theorem} (Existence of meander, co-meander and 3-d Bessel process)
\label{thm3}
For each of the following three processes $Z:= (Z_u, \le u \le 1)$ 
there is some random time $T$ such that $(B_{T+u}-B_T; 0 \leq u \leq 1)$ has the same distribution as $Z$: 
\begin{enumerate}
\item the meander $Z = (m_u; 0 \leq u \leq 1)$;
\item
the co-meander $Z = (\widetilde{m}_u; 0 \leq u \leq 1)$;
\item
the three-dimensional Bessel process $Z = (R_u; 0 \leq u \leq 1)$.
\end{enumerate}
\end{theorem}
\quad The problem of embedding Brownian bridge $b^0$ into Brownian motion is treated in a subsequent work of Pitman and Tang \cite{PTacc}. Since the proof relies heavily on {\em Palm theory} of stationary random measures, we prefer not to include it in the current work.
\begin{theorem}
\cite{PTacc}
There exists a random time $T \geq 0$ such that $(B_{T+u}-B_T; 0 \leq u \leq 1)$ has the same distribution as $(b^0_u; 0 \leq u \leq 1)$.
\end{theorem}

\quad In Question \ref{Q1}, we seek to embed a particular continuous-time process $Z$ of unit length into a Brownian path. The distribution  of $X$ resides in the infinite-dimensional space $\mathcal{C}[0,1]$ of continuous paths on $[0,1]$ starting from $0$.  So a closely related problem is whether a given subset of $\mathcal{C}[0,1]$ is hit by the path-valued moving-window process $X_t:=(B_{t+u}-B_t; 0 \leq u \leq 1)$ indexed by $t \geq 0$. We formulate this problem as follows.
\begin{question}
\label{Q2}
Given a Borel measurable subset $S \subset \mathcal{C}[0,1]$, can we find a random time $T$ such that $X_T:=(B_{T+u}-B_T; 0 \leq u \leq 1) \in S$ with probability one?
\end{question}

\quad Question \ref{Q2} involves scanning for patterns in a continuous-time process. As far as we are aware, the problem has yet been studied. Here are some examples of subsets of $\mathcal{C}[0,1]$ :
\begin{itemize}
\item
$\mathcal{E}$ is the set of excursion paths, which first returns to $0$ at time $1$, 
$$\mathcal{E}:=\{w \in \mathcal{C}[0,1]; w(t)>w(1)=0~\mbox{for}~0<t<1\};$$
\item
$\mathcal{RER}$ is the set of reflected bridge paths, i.e.
$$\mathcal{RBR}:=\{w \in \mathcal{C}[0,1]; w(t)\geq w(1)=0~\mbox{for}~0 \leq t \leq1\};$$
\item
$\mathcal{M}$ is the set of positive paths, i.e.
$$\mathcal{M}:=\{w \in \mathcal{C}[0,1]; w(t)>0~\mbox{for}~0<t \leq 1\};$$
\item
$\mathcal{BR}^{\lambda}$ is the set of bridge paths, which end at $\lambda \in \mathbb{R}$, i.e.
$$\mathcal{BR}^{\lambda}:=\{w \in \mathcal{C}[0,1]; w(1)=\lambda\};$$
\item
$\mathcal{FP}^{\lambda}$ is the set of first passage bridge paths at fixed level $\lambda<0$, i.e.
$$\mathcal{FP}^{\lambda}:=\{w \in \mathcal{C}[0,1]; w(t)>w(1)=\lambda~\mbox{for}~0 \leq t <1\};$$
\item
$\mathcal{VB}^{-}$ is the set paths of Vervaat transform of Brownian motion with a floating negative endpoint, i.e.
$$\mathcal{VB}^{-}:=\{w \in \mathcal{C}[0,1]; w(t)>w(1)~\mbox{for}~0 \leq t <1~\mbox{and}~\inf\{t > 0; w(t)<0\}>0\};$$
\item
$\mathcal{VB}^{\lambda}$ is the set of Vervaat bridge paths ending at fixed level $\lambda<0$, i.e.
$$\mathcal{VB}^{\lambda}:=\{w \in \mathcal{FP}^{\lambda}; \inf \{t> 0; w(t)<0\}>0\}=\{w \in \mathcal{VB}^{-};w(1)=\lambda\}.$$
\end{itemize}
\quad An elementary argument shows that the probability that a given set $S \subset \mathcal{C}[0,1]$ is hit by the path-valued process generated by Brownian motion is either $0$ or $1$, i.e. $$\mathbb{P}[\exists T \geq 0~\mbox{such that}~(B_{T+u}-B_T; 0 \leq u \leq 1) \in S] \in \{0,1\} \quad \mbox{for}~ S \subset \mathcal{C}[0,1].$$
Using various stochastic analysis tools, we are able to show that
\begin{theorem}(Non-existence of excursion, reflected bridge, Vervaat transform of Brownian motion, first passage bridge and Vervaat bridge paths)
\label{thm2}
\begin{enumerate}
\item\label{thm151}
Almost surely,  there is no random time $T$ such that $$(B_{T+u}-B_T; 0 \leq u \leq 1) \in \mathcal{E};$$
\item
\label{thm152}
Almost surely,  there is no random time $T$ such that $$(B_{T+u}-B_T; 0 \leq u \leq 1) \in \mathcal{RBR};$$
\item
\label{thm153}
Almost surely,  there is no random time $T$ such that $$(B_{T+u}-B_T; 0 \leq u \leq 1) \in \mathcal{VB}^{-};$$
\item
\label{thm154}
For each $\lambda<0$, almost surely,  there is no random time $T$ such that $$(B_{T+u}-B_T; 0 \leq u \leq 1) \in \mathcal{FP}^{\lambda};$$
\item
\label{thm155}
For each $\lambda<0$, almost surely,  there is no random time $T$ such that $$(B_{T+u}-B_T; 0 \leq u \leq 1) \in \mathcal{VB}^{\lambda}.$$
\end{enumerate}
\end{theorem}
\quad In the current work, we restrict ourselves to continuous paths in linear Brownian motion. However, the problem is also worth considering in the multi-dimensional case, see e.g. Section \ref{4} for further discussions. 

\quad At first glance, neither Question \ref{Q1} nor Question \ref{Q2} seems to be tractable. To gain some intuition, we start by studying the analogous problem in the random walk setting. We deal with simple symmetric random walks $SW(n)$ of length $n$ with increments $\pm 1$ starting at $0$. A typical question is how long it would take, in random walks, to observe certain element in a collection of patterns of length $n$ satisfying some common properties. More precisely, 
\begin{question}
\label{Q3}
Given for each $n$ a collection $\mathcal{A}^n$ of patterns of length $n$, what is the order of the expected waiting time $\mathbb{E}T(\mathcal{A}^n)$ until some element of $\mathcal{A}^n$ is observed in a random walk?
\end{question}
\quad  We are not aware of any previous study on pattern problems in which some natural definition of the collection of patterns is made for each $n \in \mathbb{N}$. Nevertheless, this question fits into the general theory of runs and patterns in a sequence of discrete trials. This theory dates back to work in $1940$s by Wald and Wolfowitz \cite{WW} and Mood \cite{Mood}. Since then, the subject has become important in various areas of science, including industrial engineering, biology, economics and statistics. In $1960$s, Feller \cite{Feller} treated the problem probabilistically by identifying the occurrence of a single pattern as a renewal event.  By the generating function method, the law of the occurrence times of such  patterns is entirely characterized. More advanced study of the occurrence of patterns in a collection developed in $1980$s by two different methods. Guibas and Odlyzko \cite{GO}, and Breen et al \cite{BWZ} followed the steps of Feller \cite{Feller} by studying the generating functions in pattern-overlapping regimes. An alternative approach was adopted by Li \cite{Li}, and Gerber and Li \cite{GL} using martingale arguments. We also refer readers to the book of Fu and Lou \cite{FouLou} for Markov chain embedding approach regarding multi-state trials.

\quad Techniques from the theory of patterns in i.i.d. sequences provide general strategies to Question \ref{Q3}. Here we focus on some special cases where the order of the expected waiting time is computable. As we will see later, such asymptotics help us predict the existence or non-existence of some particular patterns in Brownian motion. We introduce some discrete patterns as analogs of sets of continuous paths as above. For $\lambda \in \mathbb{R}$, let $\lambda_n=[\lambda \sqrt{n}]$ if $[\lambda \sqrt{n}]$ and $n$ have the same parity, and $\lambda_n=[\lambda \sqrt{n}]+1$ otherwise. 
We make the following definitions:
\begin{itemize}
\item
$\mathcal{E}^{2n}$ is the set of discrete positive excursions of length $2n$, whose first return to $0$ occurs at time $2n$, i.e.
$$\mathcal{E}^{2n}:=\{w \in SW(2n); w(i) > 0~\mbox{for}~1 \leq i \leq 2n-1~\mbox{and}~w(2n)=0\};$$
\item
$\mathcal{M}^{2n+1}$ is the set of positive walks of length $2n+1$, i.e.
$$\mathcal{M}^{2n+1}:=\{w \in SW(2n+1); w(i)> 0~\mbox{for}~1 \leq i \leq 2n+1\};$$
\item
$\mathcal{BR}^{\lambda,n}$ is the set of discrete bridges of length $n$, which end at $\lambda_n$ for $\lambda \in \mathbb{R}$, i.e.
$$\mathcal{BR}^{\lambda,n}:=\{w \in SW(n); w(n)=\lambda_n\};$$
\item
$\mathcal{FP}^{\lambda, n}$ for $\lambda<0$ is the set of negative first passage walks of length $n$, ending at $\lambda_n$, 
$$\mathcal{FP}^{\lambda,n}:=\{w \in SW(n); w(i) > w(n)=\lambda_n~\mbox{for}~0 \leq i \leq n-1\}.$$
\end{itemize}
\quad As in Question \ref{Q3}, $T(\mathcal{A}^n)$ is the waiting time until some pattern in $\mathcal{A}^n$ appears in the simple random walk for $\mathcal{A}^n \in \{\mathcal{E}^{2n}, \mathcal{M}^{2n+1}, \mathcal{BR}^{\lambda,n}, \mathcal{FP}^{\lambda,n}\}$. The following result answers Question \ref{Q3} in these particular cases.
\begin{theorem}
\label{thm1}
~
\begin{enumerate}
\item
\label{thm141}
There exists $C_{\mathcal{E}}>0$ such that
\begin{equation}
\label{exest}
\mathbb{E}T(\mathcal{E}^{2n}) \sim C_{\mathcal{E}} n^{\frac{3}{2}};
\end{equation}
\item
\label{thm142}
There exists $C_{\mathcal{M}}>0$ such that
\begin{equation}
\label{posest}
\mathbb{E}T(\mathcal{M}^{2n+1}) \sim C_{\mathcal{M}} n ;
\end{equation}
\item
\label{thm143}
There exist $c_{\mathcal{BR}}^{\lambda}$ and $C_{\mathcal{BR}}^{\lambda}>0$ such that
\begin{equation}
\label{brest}
c_{\mathcal{BR}}^{\lambda} n \leq \mathbb{E}T(\mathcal{BR}^{\lambda,n}) \leq  C_{\mathcal{BR}}^{\lambda} n;
\end{equation}
\item
\label{thm144}
For $\lambda<0$, there exist $c_{\mathcal{FP}}^{\lambda}$ and $C_{\mathcal{FP}}^{\lambda}>0$ such that
\begin{equation}
\label{verest}
c_{\mathcal{FP}}^{\lambda} n \leq \mathbb{E}T(\mathcal{FP}^{\lambda,n}) \leq C_{\mathcal{FP}}^{\lambda} n^{\frac{5}{4}}.
\end{equation}
\end{enumerate} 
\end{theorem}
\quad Now we explain how the asymptotics obtained in Theorem \ref{thm1} help to answer, at least formally, Question \ref{Q1} and Question \ref{Q2} in some cases. Formula
\eqref{exest}

tells that it would take on average $n^{\frac{3}{2}} \gg n$ steps to observe an excursion in simple random walk. In view of the usual scaling of random walks to converge to Brownian motion, the time scale appears to be too large. Thus we should not expect to find excursion pattern $\mathcal{E}$ in a Brownian path. However, 
in \eqref{posest} and \eqref{brest}, the typical waiting time to observe a positive walk or a bridge has the same order $n$ involved in the time scaling for convergence in distribution to Brownian motion. So we can anticipate to observe the positive-path pattern $\mathcal{M}$ and the bridge pattern $\mathcal{BR}^{\lambda}$ in Brownian motion. Finally in \eqref{verest}, there is an exponent gap in evaluating the expected waiting time to have first passage walks ending at $\lambda_n \sim [\lambda \sqrt{n}]$ for $\lambda<0$. In this case, we do not know whether it would take asymptotically $n$ steps or much longer to first observe such patterns. This prevents us from predicting the existence of the first passage bridge pattern $\mathcal{FP}^{\lambda}$ in Brownian motion.

\quad The scaling arguments used in the last paragraph are quite intuitive but not rigorous since we are not aware of any theory which would justify the existence or non-existence of continuous paths by taking  limits from the discrete setting.\\\\
\textbf{Organization of the paper:} The rest of the paper is organized as follows.
\begin{itemize}
\item 
In Section \ref{2}, we study the asymptotic behavior of the expected waiting time for discrete patterns. There Theorem \ref{thm1} is proved.
\item
In Section \ref{3}, we are devoted to the analysis of continuous paths/processes in Brownian motion. Proofs of Theorem \ref{thm3} and Theorem \ref{thm2} are provided.
\item
In Section \ref{4}, we discuss the potential theory for continuous paths in Brownian motion. 
\end{itemize}
A selection of open problems is presented in Section \ref{25} and Section \ref{4}.
\section{Expected waiting time for discrete patterns}
\label{2}
\quad The current section is devoted to the proof of Theorem \ref{thm1}. We study the expected waiting time for some collection of patterns $\mathcal{A}^n \in \{\mathcal{E}^{2n}, \mathcal{M}^{2n+1},  \mathcal{BR}^{\lambda,n},\mathcal{FP}^{\lambda,n}\}$ as defined in the introduction. Note that all patterns in the collection $\mathcal{A}^n$ have the same length. More precisely, we are interested in the asymptotic behavior of $\mathbb{E}T(\mathcal{A}^n)$ as the common length of patterns $n \rightarrow \infty$.

\quad We start by recalling the general strategy to compute the expected waiting time for discrete patterns in simple random walks. For $\mathcal{A}^n:=\{A_{1}^n, \cdots, A_{k(\mathcal{A}^n)}^n\}$, denote $T(A_{i}^n)$ the waiting time for the $A_{i}^n$ renewal and $T(\mathcal{A}^n)$ the waiting time till some $A_{i}^n \in \mathcal{A}^n$ is renewed, see e.g. Feller \cite[Chapter XIII]{Feller} for general background on renewal strings/patterns.

\quad Define the matching matrix $M(\mathcal{A}^n)$, which accounts for the overlapping phenomenon among patterns within the collection $\mathcal{A}^n$. The coefficients are given by
\begin{equation}
\label{match}
M(\mathcal{A}^n)_{ij}:=\sum_{l=0}^{n-1}\frac{\epsilon_l(A_i^n,A_j^n)}{2^l} \quad \mbox{for}~ 1 \leq i,j \leq k(\mathcal{A}^n),
\end{equation}
where $\epsilon_{l}(A_i^n,A_j^n)$ is defined for $A_i^n:=A_{i1}^n \cdots A_{in}^n$ and $A_j^n:=A_{j1}^n \cdots A_{jn}^n$ as
\begin{equation}
\label{eps}
\epsilon_l(A_i^n,A_j^n) := \left\{ \begin{array}{rcl}
         1 & \quad \mbox{if}~
         A_{i1}^n=A_{j1+l}^n , \cdots , A_{in-l}^n=A_{jn}^n \\
         0 & \mbox{otherwise},
                \end{array}\right.
\end{equation}
for $0 \leq l \leq n-1$. Note that in general for $i \neq j$, $M(\mathcal{A}^n)_{ij} \neq M(\mathcal{A}^n)_{ji}$ and hence the matching matrix $M(\mathcal{A}^n)$ is not necessarily symmetric. The following result, which can be read from Breen et al \cite{BWZ} is the main tool to study the expected waiting time for the collection of patterns.
\begin{theorem} \cite{BWZ}
\label{BWZ}
\begin{enumerate}
\item
The matching matrix $M(\mathcal{A}^n)$ is invertible and the expected waiting times for patterns in $\mathcal{A}^n:=\{A_{1}^n, \cdots, A_{k(\mathcal{A}^n)}^n\}$ are given by
\begin{equation}
\label{ind}
\left(\frac{1}{\mathbb{E}T(A_{1}^n)}, \cdots , \frac{1}{\mathbb{E}T(A_{k(\mathcal{A}^n)}^n)} \right)^T = \frac{1}{2^n} M(\mathcal{A}^n)^{-1}  \left(1, \cdots , 1 \right)^T;
\end{equation}
\item
The expected waiting time till one of the patterns in $\mathcal{A}^n$ is observed is given by
\begin{equation}
\label{total}
\frac{1}{\mathbb{E}T(\mathcal{A}^n)} = \sum_{l=1}^{k(\mathcal{A}^n)} \frac{1}{\mathbb{E}T(A_{l}^n)}=\frac{1}{2^n} (1, \cdots , 1) M(\mathcal{A}^n)^{-1} (1, \cdots, 1)^T.
\end{equation}
\end{enumerate}
\end{theorem}
\quad In Section \ref{21}, we apply the previous theorem to obtain the expected waiting time for discrete excursions $\mathcal{E}^{2n}$, i.e. \eqref{thm141} of Theorem \ref{thm1}. The same problem for positive walks $\mathcal{M}^{2n+1}$, bridge paths $\mathcal{BR}^{0,2n}$ and first passage walks $\mathcal{FP}^{\lambda,n}$ through $\lambda_n \sim \lambda \sqrt{n}$, i.e. \eqref{thm142}$-$\eqref{thm144} of Theorem \ref{thm1}, is studied respectively in Section \ref{22}$-$\ref{24}. Finally, we discuss the exponent problem for certain discrete patterns in Section \ref{25}. 
\subsection{Expected waiting time for discrete excursions}
\label{21}
For $n \in \mathbb{N}$, the number of discrete excursions of length $2n$ is equal to the $n-1^{th}$ Catalan number, see e.g. Stanley \cite[Exercise $6.19$ (i)]{Stanley}, i.e.
\begin{equation}
\label{catalan}
k(\mathcal{E}^{2n}):=\# \mathcal{E}^{2n} =\frac{1}{n} \binom{2n-2}{n-1} \sim \frac{1}{4 \sqrt{\pi}} 2^{2n} n^{-\frac{3}{2}}.
\end{equation}
\quad Note that discrete excursions never overlap since the starting point and the endpoint are the only two minima. We have then $\epsilon(E_{i}^n,E_{j}^n)=\delta_{ij}$ for $1 \leq i, j \leq k(\mathcal{E}^{2n})$ by \eqref{eps}. Thus, the matching matrix defined as in \eqref{match} for discrete excursions $\mathcal{E}^{2n}$ has the simple form
\begin{equation*}
M(\mathcal{E}^{2n})=I_{k(\mathcal{E}^{2n})} \quad (k(\mathcal{E}^{2n}) \times k(\mathcal{E}^{2n})~\mbox{identity matrix}).
\end{equation*} 
According to Theorem \ref{BWZ},
\begin{equation}
\label{exest2}
\forall 1 \leq i \leq k(\mathcal{E}^{2n}),~\mathbb{E}T(E_{i}^n)=\frac{1}{2^{2n}} \quad \mbox{and} \quad \mathbb{E}T(\mathcal{E}^n) =\frac{2^{2n}}{k(\mathcal{E}^{2n})} \sim 4 \sqrt{\pi} n^{\frac{3}{2}},
\end{equation}
where $k(\mathcal{E}^{2n})$ is given as in \eqref{catalan}. Then \eqref{exest} follows from \eqref{exest2} with $C_{\mathcal{E}}=4 \sqrt{\pi}$.  $\square$
\subsection{Expected waiting time for positive walks}
\label{22}
Let $n \in \mathbb{N}$. It is well-known that the number of non-negative walks of length $2n+1$ is $\binom{2n}{n}$, see e.g. Larbarbe and Marckert \cite{LM} and Leeuwen \cite{Leeuwen} for modern proofs. Thus the number of positive walks of length $2n+1$ is given by
\begin{equation}
\label{enumpositive}
k(\mathcal{M}^{2n+1}):=\# \mathcal{M}^{2n+1}= \binom{2n}{n} \sim \frac{1}{\sqrt{\pi}} 2^{2n}n^{-\frac{1}{2}}.
\end{equation}
\quad Now consider the matching matrix $M(\mathcal{M}^{2n+1})$ defined as in \eqref{match} for positive walks $\mathcal{M}^{2n+1}$. $M(\mathcal{M}^{2n+1})$ is no longer diagonal since there are overlaps among positive walks. Nevertheless, the following result unveils the particular structure of such matrix.
\begin{lemma}
\label{stochmatrix}
$M(\mathcal{M}^{2n+1})$ is the multiple of some right stochastic matrix (whose row sums are equal to $1$) with multiplicity
\begin{equation}
\label{multiplicity}
1+\sum_{l=1}^{2n} \frac{k(\mathcal{M}^l)}{2^l} \sim \frac{2}{\sqrt{\pi}} \sqrt{n},
\end{equation}
where $k(\mathcal{M}^l):=\# \mathcal{M}^l$ is the number of positive walks of length $l$ as in \eqref{enumpositive}.
\end{lemma}
\textbf{Proof:} Let $1 \leq i \leq k(\mathcal{M}^{2n+1})$ and consider the sum of the $i^{th}$ row 
\begin{align}
\label{fubini}
\sum_{j=1}^{k(\mathcal{M}^{2n+1})} M(\mathcal{M}^{2n+1})_{ij} :&= \sum_{j=1}^{k(\mathcal{M}^{2n+1})} \sum_{l=0}^{2n} \frac{\epsilon_l(M_i^{2n+1},M_j^{2n+1})}{2^l} \notag\\
&= \sum_{l=0}^{2n} \frac{1}{2^l} \sum_{j=1}^{k(\mathcal{M}^{2n+1})} \epsilon_l(M_i^{2n+1},M_j^{2n+1}), 
\end{align}
where for $0 \leq l \leq 2n$ and $M_i^{2n+1}, M_j^{2n+1} \in \mathcal{M}^{2n+1}$, $\epsilon_l(M_i^{2n+1},M_j^{2n+1})$ is defined as in \eqref{eps}. Note that $\epsilon_{0}(M_i^{2n+1},M_j^{2n+1})=1$ if and only if $i=j$. Thus,
\begin{equation}
\label{0case}
\sum_{j=1}^{k(\mathcal{M}^{2n+1})} \epsilon_0(M_i^{2n+1},M_j^{2n+1})=1.
\end{equation}
In addition, for $1 \leq l \leq 2n$,
\begin{multline*}
\sum_{j=1}^{k(\mathcal{M}^{2n+1})} \epsilon_l(M_i^{2n+1},M_j^{2n+1}) = \\ \#\{M_j^{2n+1} \in \mathcal{M}^{2n+1}; M_{i1}^{2n+1}=M_{j1+l}^{2n+1}, \cdots , M_{i2n+1-l}^{2n+1}=M_{j2n+1}^{2n+1}\}.
\end{multline*}
Note that given $M_{i1}^{2n+1}=M_{j1+l}^{2n+1}, \cdots , M_{i2n+1-l}^{2n+1}=M_{j2n+1}^{2n+1}$, which implies that $M_{j1+l}^{2n+1} \cdots M_{j2n+1}^{2n+1}$ is a positive walk of length $2n-l+1$, we have
$$M_j^{2n+1} \in \mathcal{M}^{2n+1} \Longleftrightarrow M_{j1}^{2n+1} \cdots M_{jl}^{2n+1}~\mbox{is a positive walk of length}~l.$$
Therefore, for $1 \leq l \leq 2n$,
\begin{equation}
\label{lcase}
\sum_{j=1}^{k(\mathcal{M}^{2n+1})} \epsilon_l(M_i^{2n+1},M_j^{2n+1}) = k(\mathcal{M}^l).
\end{equation}
In view of \eqref{fubini}, \eqref{0case} and \eqref{lcase}, we obtain for all $1 \leq i \leq k(\mathcal{M}^{2n+1})$, the sum of $i^{th}$ row of $M(\mathcal{M}^{2n+1})$ is given by \eqref{multiplicity}. Furthermore, by \eqref{enumpositive}, we know that $k(\mathcal{M}^l) \sim \frac{1}{\sqrt{2 \pi}}2^l l^{-\frac{1}{2}}$ as $l \rightarrow \infty$, which yields the asymptotics $\frac{2}{\sqrt{\pi}}\sqrt{n}$.  $\square$\\

\quad Now by $(1)$ in Theorem \ref{BWZ}, $M(\mathcal{M}^{2n+1})$ is invertible and the inverse $M(\mathcal{M}^{2n+1})^{-1}$ is as well the multiple of certain right stochastic matrix with multiplicity
$$\left(1+\sum_{l=1}^{n-1} \frac{k(\mathcal{M}^l)}{2^l}\right)^{-1} \sim \frac{\sqrt{\pi}}{2 \sqrt{n}}.$$
Then using \eqref{total}, we obtain
\begin{equation}
\label{positivewalk}
\mathbb{E}T(\mathcal{M}^{2n+1}) = \frac{2^{2n+1}}{\left(1+\sum_{l=1}^{n-1} \frac{k(\mathcal{M}^l)}{2^l}\right)^{-1}k(\mathcal{M}^{2n+1})} \sim 4n.
\end{equation}
Therefore, \eqref{posest} follows from \eqref{positivewalk} by taking $C_{\mathcal{M}}=4$.  $\square$
\subsection{Expected waiting time for bridge paths}
\label{23}
In this part, we deal with the expected waiting time for the set of discrete bridges. In order to simplify the notations, we focus on the set of bridges of length $2n$ which ends at $\lambda=0$, i.e. $\mathcal{BR}^{0,2n}$. Note that the result in the general case for $\mathcal{BR}^{\lambda,n}$, where $\lambda \in \mathbb{R}$,  can be derived in a similar way. Using Theorem \ref{BWZ}, we prove a weaker version of \eqref{brest}, i.e. there exist $\widetilde{c}_{\mathcal{BR}}^{0}$ and $C_{\mathcal{BR}}^{0}>0$ such that
\begin{equation}
\label{brestbis}
\widetilde{c}_{\mathcal{BR}}^{0} n^{\frac{1}{2}} \leq \mathbb{E}T(\mathcal{BR}^{0,n}) \leq C_{\mathcal{BR}}^{0} n.
\end{equation}
Compared to \eqref{brest}, there is an exponent gap in \eqref{brestbis} and the lower bound is not optimal. Nevertheless, the lower bound of \eqref{brest} follows a soft argument by scaling limit, Proposition \ref{scaling}. We defer the discussion to Section \ref{25}. It is standard that the number of discrete bridges of length $2n$ is 
\begin{equation}
\label{enumbr}
k(\mathcal{BR}^{0,2n}):=\# \mathcal{BR}^{0,2n}=\binom{2n}{n} \sim \frac{1}{\sqrt{\pi}} 2^{2n} n^{-\frac{1}{2}}.
\end{equation}
\quad Denote $\mathcal{BR}^{0,2n}:=\{BR^{2n}_1, \cdots, BR^{2n}_{k(\mathcal{BR}^{0,2n})}\}$ and $M(\mathcal{BR}^{0,2n})$ the matching matrix of $\mathcal{BR}^{0,2n}$. We first establish the LHS estimate of \eqref{brestbis}. According to \eqref{ind}, we have
\begin{equation}
\label{bbq}
(1, \cdots, 1) M(\mathcal{BR}^{0,2n}) \left(\frac{1}{\mathbb{E}T(BR^{2n}_1)}, \cdots, \frac{1}{\mathbb{E}T(BR^{2n}_{k(\mathcal{BR}^{0,2n})})}\right)^T = \frac{k(\mathcal{BR}^{0,2n})}{2^{2n}}.
\end{equation}
Note that the matching matrix $\mathcal{M}(\mathcal{BR}^{0,2n})$ is non-negative with diagonal elements
$$M(\mathcal{BR}^{0,2n})_{ii} \geq \epsilon_0(BR^{2n}_i, BR^{2n}_i)=1,$$
for $1 \leq i \leq k(\mathcal{BR}^{0,2n})$. As a direct consequence, the column sums of $M(\mathcal{BR}^{0,2n})$ is larger or equal to $1$. Then by \eqref{total} and \eqref{bbq},
$$\mathbb{E}T(\mathcal{BR}^{0,2n}) \geq \frac{2^{2n}}{k(\mathcal{BR}^{0,2n})} \sim \sqrt{\pi n},$$
where $k(\mathcal{BR}^{0,2n})$ is defined as in \eqref{enumbr}. Take then $\widetilde{c}_{\mathcal{BR}}^{0}=\sqrt{\pi}-\epsilon$ for arbitrary small $\epsilon>0$.

\quad Now we establish the RHS estimate of \eqref{brestbis}. In view of \eqref{bbq}, it suffices to work out an upper bound for the column sums of $M(\mathcal{BR}^{0,2n})$. Similarly as in \eqref{fubini}, for $1 \leq j \leq k(\mathcal{BR}^{0,2n})$,
\begin{equation}
\label{columnbr}
\sum_{i=1}^{k(\mathcal{BR}^{0,2n})} M(\mathcal{BR}^{0,2n})_{ij} = 1 + \sum_{l=1}^{2n-1} \frac{1}{2^l} \sum_{i=1}^{k(\mathcal{BR}^{0,2n})} \epsilon_l(BR^{2n}_i,BR^{2n}_j),
\end{equation}
and
\begin{align}
\sum_{i=1}^{k(\mathcal{BR}^{0,2n})} \epsilon_l(BR^{2n}_i,BR^{2n}_j)
 &= \# \{BR_i^{2n} \in \mathcal{BR}^{0,2n}; BR^{2n}_{i1}=BR^{2n}_{j1+l}, \cdots , BR^{2n}_{in-l}=BR^{2n}_{jn}\}. \notag\\
 &= \# \{\mbox{discrete bridges of length}~l~\mbox{which ends at}~\sum_{k=1}^{n-l} BR^{2n}_{jk}\} \notag\\
 & \label{goodest}= \binom{l}{\frac{l+\sum_{k=1}^{n-l} BR^{2n}_{jk}}{2}} \leq \binom{l}{[\frac{l}{2}]},
 \end{align}
where the last inequality is due to the fact that $\binom{l}{k} \leq \binom{l}{[l]/2}$ for $0 \leq k \leq l$. By \eqref{columnbr} and \eqref{goodest}, the column sums of $M(\mathcal{BR}^{0,2n})$ is bounded from above by $$1+ \sum_{l=0}^{2n-1} \frac{1}{2^l} \binom{l}{[\frac{l}{2}]} \sim \frac{4}{\sqrt{\pi}} n^{\frac{1}{2}}.$$

Again by \eqref{total} and \eqref{bbq},
$$\mathbb{E}T(\mathcal{BR}^{0,2n}) \leq 2^{2n} \frac{4n^{\frac{1}{2}}/\sqrt{\pi}}{k(\mathcal{BR}^{0,2n})} \sim 4 n.$$
Hence we take $C^{0}_{\mathcal{BR}}=4+\epsilon$ for arbitrary small $\epsilon>0$.  $\square$
\subsection{Expected waiting time for first passage walks}
\label{24}
We consider the expected waiting time for first passage walks through $\lambda_n \sim \lambda \sqrt{n}$ for $\lambda<0$. Following Feller \cite[Theorem $2$, Chapter III.$7$]{Feller}, the number of patterns in $\mathcal{FP}^{\lambda,n}$ is
\begin{equation}
\label{enumfpb}
k(\mathcal{FP}^{\lambda,n}):= \# \mathcal{FP}^{\lambda, n}= \frac{\lambda_n}{n} \binom{n}{\frac{n+\lambda_n}{2}} \sim \lambda \exp \left( -\frac{\lambda^2}{2}\right)\sqrt{\frac{2}{\pi}}  2^n n^{-1}.
\end{equation}
\quad For $\mathcal{FP}^{\lambda,n}:=\{FP^n_1, \cdots , FP^n_{k(\mathcal{FP}^{\lambda,n})}\}$ and $M(\mathcal{FP}^{\lambda,n})$ the matching matrix for $\mathcal{FP}^{\lambda,n}$, we have, by \eqref{ind}, that
\begin{equation}
\label{consind}
(1, \cdots, 1) M(\mathcal{FP}^{\lambda,n}) \left(\frac{1}{\mathbb{E}T(FP^n_1)}, \cdots, \frac{1}{\mathbb{E}T(FP^n_{k(\mathcal{FP}^{\lambda,n})})}\right)^T = \frac{k(\mathcal{FP}^{\lambda,n})}{2^n}.
\end{equation}
The LHS bound of \eqref{verest} can be derived in a similar way as in Section \ref{23}, i.e.
\begin{equation*}
\label{leftest}
\mathbb{E}T(\mathcal{FP}^{\lambda,n}) \geq \frac{2^n}{k(\mathcal{FP}^{\lambda,n})} \sim \sqrt{\frac{\pi}{2 \lambda^2}}\exp \left( \frac{\lambda^2}{2}\right) n, 
\end{equation*}
where $k(\mathcal{FP}^{\lambda,n})$ is defined as in \eqref{enumfpb}. We take $c^{\lambda}_{\mathcal{FP}}=\sqrt{\frac{\pi}{2 \lambda^2}}\exp(\frac{\lambda^2}{2})-\epsilon$ for arbitrary small $\epsilon>0$.

\quad For the upper bound of \eqref{verest}, we aim to obtain an upper bound for the column sums of $M(\mathcal{FP}^{\lambda,n})$. Note that for $1 \leq j \leq k_{\mathcal{FP}^n}^{\lambda}$,
\begin{equation}
\label{columnsum}
\sum_{i=1}^{k(\mathcal{FP}^{\lambda,n})} M(\mathcal{FP}^{\lambda,n})_{ij} = 1 + \sum_{l=1}^{n-1} \frac{1}{2^l} \sum_{i=1}^{k(\mathcal{FP}^{\lambda,n})} \epsilon_l(FP_i,FP_j)
\end{equation}
and
$$\sum_{i=1}^{k(\mathcal{FP}^{\lambda,n})} \epsilon_l(FP_i^n,FP_j^n) = \# \{FP_i^n \in \mathcal{FP}^{\lambda,n}; FP_{i1}^n=FP_{j1+l}^n, \cdots , FP_{in-l}^n=FP_{jn}^n\}.$$
\quad Observe that $\{FP_i^n \in \mathcal{FP}^{\lambda,n}; FP_{i1}^n=FP_{j1+l}^n, \cdots , FP_{in-l}^n=FP_{jn}^n\} \neq \emptyset$ if and only if $\sum_{k=1}^{l}FP_{jk}^n<0$ (otherwise $\sum_{k=1}^{n-l}FP_{ik}^n=\sum_{k=1+l}^nFP_{jk}^n=\lambda_n - \sum_{k=1}^l FP_{jk}^n < \lambda_n$, which implies $FP_i^n \notin \mathcal{FP}^{\lambda,n}$). Then given $FP_{i1}^n=FP_{j1+l}^n, \cdots , FP_{in-l}^n=FP_{jn}^n$ and $\sum_{k=1}^{l}FP_{jk}^n<0$, 
\begin{multline*}
FP_i^n \in \mathcal{FP}^{\lambda,n} \Longleftrightarrow \\ FP_{in-l+1}^n\cdots FP_{in}^n~\mbox{is a first passage walk of length}~l~\mbox{through}~\sum_{k=1}^{l}FP_{jk}^n<0.
\end{multline*} 
Therefore, for $1 \leq l \leq n-1$ and $1 \leq j \leq k_{\mathcal{FP}^{n}}^{\lambda}$,
\begin{equation}
\label{indicator}
\sum_{i=1}^{k(\mathcal{FP}^{\lambda,n})} \epsilon_l(FP_i^n,FP_j^n) = 1_{\sum_{k=1}^l FP_{jk}^n<0} \frac{|\sum_{k=1}^l FP_{jk}^n|}{l} \binom{l}{\frac{l+\sum_{k=1}^l FP_{jk}^n}{2}}.
\end{equation}
\quad From the above discussion, it is easy to see for $1 \leq j \leq k(\mathcal{FP}^{\lambda,n})$,
$$\sum_{i=1}^{k(\mathcal{FP}^{\lambda,n})} M(\mathcal{FP}^{\lambda,n})_{ij} \leq \sum_{i=1}^{k(\mathcal{FP}^{\lambda,n})} M(\mathcal{FP}^{\lambda,n})_{ij^{*}},$$
where $FP_{j^{*}}^n$ is defined as follows: $FP_{j^{*}k}^n = -1$ if $1 \leq k \leq \lambda_n-1$; $\lambda_n-1<k \leq n-1$ and $k-\lambda_n$ is odd; $k=n$. Otherwise $FP_{j^{*}k}^n = 1$.
\begin{figure}[h]
\includegraphics[width=0.6 \textwidth]{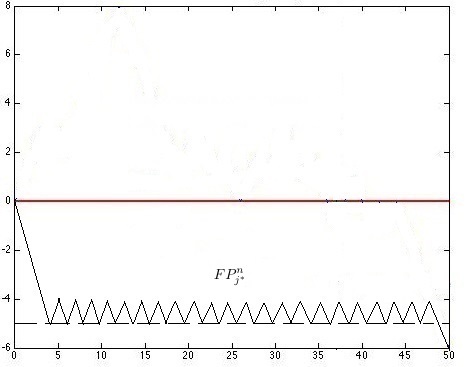}
\caption{Extreme patterns $FP_{j^{*}}^n$.}
\end{figure}\\
The rest of this part is devoted to estimating $\sum_{i=1}^{k(\mathcal{FP}^{\lambda,n})} M(\mathcal{FP}^{\lambda,n})_{ij^{*}}$. By \eqref{columnsum} and \eqref{indicator},
\begin{align}
\label{FPlarge}
\sum_{i=1}^{k(\mathcal{FP}^{\lambda,n})} M(\mathcal{FP}^{\lambda,n})_{ij^{*}}&=\sum_{l=0}^{|\lambda_n|-1}\frac{1}{2^l}+\underset{l-|\lambda_n|~\mbox{\scriptsize{odd}}}{\sum_{l=\lambda_n}^{n-1}}  \frac{|\lambda_n|-1}{l \cdot 2^l} \binom{l}{\frac{l-|\lambda_n|+1}{2}}+\underset{l-|\lambda_n|~\mbox{\scriptsize{even}}}{\sum_{l=\lambda_n}^{n-1}} \frac{|\lambda_n|-2}{l \cdot 2^l} \binom{l}{\frac{l-|\lambda_n|+2}{2}} \notag\\
& \leq 2+ |\lambda_n|  \sum_{l=|\lambda_n|}^{n-1} \frac{1}{2^ll} \binom{l}{[\frac{l}{2}]}  \sim \sqrt{\frac{8 \lambda}{\pi}} n^{\frac{1}{4}}. \notag
\end{align}
Thus, the column sums of $M(\mathcal{FP}^{\lambda,n})$ is bounded from above by $\sqrt{\frac{8 \lambda}{\pi}} n^{\frac{1}{4}}$. By \eqref{total} and \eqref{consind}, 
\begin{equation*}
\label{rightest}
\mathbb{E}T(\mathcal{FP}^{\lambda,n}) \leq \frac{2^n \sqrt{8 \lambda/\pi}n^{\frac{1}{4}}}{k(\mathcal{FP}^{\lambda,n})} \sim \sqrt{\frac{4}{\lambda}} \exp \left(\frac{\lambda^2}{2}\right) n^{\frac{5}{4}}. 
\end{equation*}
We take then $C_{\mathcal{FP}}^{\lambda}=\sqrt{\frac{4}{\lambda}} \exp \left(\frac{\lambda^2}{2}\right)+\epsilon$ for arbitrary small $\epsilon>0$.  $\square$
\subsection{Exponent gaps for $\mathcal{BR}^{\lambda,n}$ and $\mathcal{FP}^{\lambda,n}$}
\label{25}
It can be inferred from \eqref{brestbis} (resp. \eqref{verest}) that the expected waiting time for $\mathcal{BR}^{\lambda,n}$ where $\lambda \in \mathbb{R}$ (resp. $\mathcal{FP}^{\lambda,n}$ where $\lambda<0$) is bounded from below by order $n^{\frac{1}{2}}$ (resp. $n$) and from above by order $n$ (resp. $n^{\frac{5}{4}}$). The exponent gap in the estimates of first passage walks $\mathcal{FP}^{\lambda,n}$ is frustrating, since we do not know whether the waiting time is exactly of order $n$, or is of order $\gg n$. This prevents the prediction of the existence of first passage bridge patterns $\mathcal{FP}^{\lambda}$ in Brownian motion. 

\quad From \eqref{total}, we see that the most precise way to compute $\mathbb{E}T(\mathcal{BR}^{\lambda,n})$ and $\mathbb{E}T(\mathcal{FP}^{\lambda,n})$ consists in evaluating the sum of all entries in the inverse matching matrices $M(\mathcal{BR}^{\lambda,n})^{-1}$ and $M(\mathcal{FP}^{\lambda,n})^{-1}$. But the task is difficult since the structures of $M(\mathcal{BR}^{\lambda,n})$ and $M(\mathcal{FP}^{\lambda,n})$ are complex, at least more inoperable than those as of $M({\mathcal{E}^{2n}})$ and $M(\mathcal{M}^{2n+1})$. We do not understand well the exact form of its inverse matrix $M(\mathcal{BR}^{\lambda,n})^{-1}$ and $M(\mathcal{FP}^{\lambda,n})^{-1}$.

\quad The technique used in Section \ref{23} and Section \ref{24} is to bound the column sums of the matching matrix $M(\mathcal{BR}^{\lambda,n})$ (resp. $M(\mathcal{FP}^{\lambda,n})$). More precisely, we have proved that 
\begin{equation}
\label{roughest1}
\mathcal{O}(1) \leq~\mbox{column sums of}~M(\mathcal{BR}^{\lambda,n}) \leq  \mathcal{O}(n^{\frac{1}{2}}) \quad \mbox{for each fixed}~\lambda \in \mathbb{R};
\end{equation}
\begin{equation}
\label{roughest2}
\mathcal{O}(1) \leq~\mbox{column sums of}~M(\mathcal{FP}^{\lambda,n}) \leq  \mathcal{O}(n^{\frac{1}{4}})\quad \mbox{for each fixed}~\lambda <0.
\end{equation}
\quad For the bridge pattern $\mathcal{BR}^{0,2n}$, the LHS bound of \eqref{roughest1} is obtained by any excursion path of length $2n$, while the RHS bound of \eqref{roughest1} is achieved by the sawtooth path, with consecutive $\pm 1$ increments. In the first passage pattern $\mathcal{FP}^{\lambda,n}$ where $\lambda<0$, the LHS bound of \eqref{roughest2} is attained by some excursion-like paths, which start by excursions and go linearly to $\lambda \sqrt{n}<0$ at the end. The RHS bound of \eqref{roughest2} is given by the extreme pattern defined in Section \ref{24}, see e.g. Figure $1$.

\quad However, all above estimations are not accurate, since there are only few columns in $\mathcal{BR}^{\lambda,n}$ (resp. $\mathcal{FP}^{\lambda,n}$) which sum up either to $\mathcal{O}(1)$ (resp. $\mathcal{O}(1)$) or to $\mathcal{O}(n^{\frac{1}{2}})$ (resp. $\mathcal{O}(n^{\frac{1}{4}})$). Thus,
\begin{openpb} Using the matching matrix method as in Theorem \ref{BWZ},
\begin{enumerate}
\item
determine the exact asymptotics for $\mathbb{E}T(\mathcal{BR}^{\lambda,n})$ where $\lambda \in \mathbb{R}$, as $n \rightarrow \infty$;
\item
determine the exact asymptotics for $\mathbb{E}T(\mathcal{FP}^{\lambda,n})$ where $\lambda<0$, as $n \rightarrow \infty$.
\end{enumerate}
\end{openpb}

\quad As we prove below, for $\lambda \in \mathbb{R}$, $\mathbb{E}T(\mathcal{BR}^{\lambda,n}) \asymp n$ by a scaling limit argument. Nevertheless, to obtain this result only by discrete analysis would be of independent interest. The following table provides the simulations of the expected waiting time $\mathbb{E}T(\mathcal{FP}^{-1,n})$ for some large $n$. The result strongly suggests the existence of a real number $\zeta \in (1,\frac{5}{4})$ such that $\mathbb{E}T(\mathcal{FP}^{-1,n}) \asymp n^{\zeta}$. This is confirmed by Theorem \ref{thm2} \eqref{thm154}, that is we cannot find a first passage bridge with fixed negative endpoint in Brownian motion.
\begin{center}
    \begin{tabular}{| c | c | c | c | c | c|}
    \hline
    $n$ & $100$ & $200$ & $500$ & $1000$ & $2000$ \\ \hline
    $\mathbb{E}T(\mathcal{FP}_{-1}^n)$ & $179.8050$ & $358.2490$ & $893.0410$ & $1800.0020$ & $3682.0220$  \\ \hline
    Estimated $\zeta$&  & $0.9945$ & $0.9968$ & $1.0112$ & $1.0375$  \\ \hline
    \end{tabular}\\
    \end{center}
{TABLE $1$. Estimation of $\zeta$ by $\log \frac{\mathbb{E}T(\mathcal{FP}_{-1}^{n_2})}{\mathbb{E}T(\mathcal{FP}_{-1}^{n_1})}/ \log (\frac{n_2}{n_1}$), where $n_2$ is the next to $n_1$ in the table.}

\quad Now let us focus on the lower bound \eqref{brest} of expected waiting time for bridge pattern $\mathcal{BR}^{0}$. For $n \in 2 \mathbb{N}$, we run a simple random walk $(RW_k)_{k \in \mathbb{N}}$ until the first level bridge of length $n$ appears. That is, we consider
\begin{equation}
\label{disbl}
(RW_{F_n+k}-RW_{F_n})_{0 \leq k \leq n}, \quad \mbox{where}~F_n:=\inf\{k \geq 0; RW_{k+n}=RW_k\}.
\end{equation}
For simplicity, let $RW_k$ for non-integer $k$ be defined by the usual linear interpolation of a simple random walk. For background on the weak convergence in $\mathcal{C}[0,1]$, we refer readers to Billingsley \cite[Chapter 2]{Bill}.
\begin{proposition}
\label{scaling}
$$\left(\frac{RW_{F_n+nu}-RW_{F_n}}{\sqrt{n}}; 0 \leq u \leq 1 \right)~\mbox{converges weakly in}~\mathcal{C}[0,1]~\mbox{to the bridge-like process}$$
\begin{equation}
\label{bridgelike}
(B_{F+u}-B_F; 0 \leq u \leq 1), \quad \mbox{where}~F:=\inf\{t>0; B_{t+1}-B_{t}=0\}.
\end{equation}
\end{proposition}
\quad The process $(S_t:=B_{t+1}-B_{t}; t \geq 0)$ is a stationary Gaussian process, first studied by Slepian \cite{Slepian} and Shepp \cite{Shepp}. The following result, which can be found in Pitman and Tang \cite[Lemma $2.3$]{PTacc}, is needed for the proof of Proposition \ref{scaling}. 
\begin{lemma} \cite{SheppGaussian,PTacc}
\label{abszero}
For each fixed $t \geq 0$, the distribution of $(S_u; t \leq u \leq t+1)$ is mutually absolutely continuous with respect to the distribution of
\begin{equation}
\label{Btil}
(\widetilde{B}_u:=\sqrt{2}(\xi+B_u); t \leq u \leq t+1),
\end{equation}
where $\xi \sim \mathcal{N}(0,1)$. In particular, the distribution of the Slepian zero set restricted to $[t,t+1]$, i.e. $\{u \in [t,t+1]; S_u=0\}$ is mutually absolutely continuous with respect to that of $\{u \in [t,t+1]; \xi+B_u=0\}$, the zero set of Brownian motion starting at $\xi \sim \mathcal{N}(0,1)$.
\end{lemma}
\textbf{Proof of Proposition \ref{scaling}:} Let $\mathbb{P}^{\bf W}$ be Wiener measure on $\mathcal{C}[0,\infty)$. Let $\mathbb{P}^{{\bf S}}$ (resp. $\mathbb{P}^{\bf \widetilde{{\bf W}}}$) be the distribution of the Slepian process $S$ (resp. the distribution of $\widetilde{B}$ defined as in \eqref{Btil}). We claim that $$F:=\inf\{t \geq 0; w_{t+1}=w_t\},$$ is a  functional of the coordinate process $w:=\{w_t; t \geq 0\} \in \mathcal{C}[0,\infty)$ that is continuous $\mathbb{P}^{\bf W}$ a.s. Note that the distribution of $(x_t:=w_{t+1}-w_{t}; t \geq 0)$ under $\mathbb{P}^{\bf W}$ is the same as that of $(w_t;t \geq 0)$ under $\mathbb{P}^{\bf S}$. In addition, $x \in \mathcal{C}[0,\infty)$ is a functional of $w \in \mathcal{C}[0,\infty)$ that is continuous $\mathbb{P}^{\bf W}$ a.s. By composition, it is equivalent to show that
$$F':=\inf\{t \geq 0; w_t=0\},$$ is a functional of $w \in \mathcal{C}[0,\infty)$ that is continuous $\mathbb{P}^{\bf S}$ a.s. Consider the set
$$\mathcal{Z}:=\{F'~\mbox{is not a continuous functional of}~w \in \mathcal{C}[0,\infty)\} = \cup_{p \in \mathbb{Q}} \mathcal{Z}_p,$$
where $\mathcal{Z}_p:=\{G' \in [p,p+1]~\mbox{and}~F'~\mbox{is not a continuous functional of}~w \in \mathcal{C}[0,\infty)\}.$ It is obvious that $\mathbb{P}^{\widetilde{{\bf W}}}(\mathcal{Z})=0$ and thus $\mathbb{P}^{\widetilde{{\bf W}}}(\mathcal{Z}_p)=0$ for all $p \geq 0$. By Lemma \ref{abszero}, $\mathbb{P}^{\bf S}$ is locally absolutely continuous relative to $\mathbb{P}^{\widetilde{{\bf W}}}$, which implies that $\mathbb{P}^{\bf S}(\mathcal{Z}_p)=0$ for all $p \geq 0$. As a countable union of null events, $\mathbb{P}^{\bf S}(\mathcal{Z})=0$, and the claim is proved. Thus, the mapping
$$\Xi_F: \mathcal{C}[0,\infty) \ni (w_t; t \geq 0) \longrightarrow (w_{F+u}-w_F; 0 \leq u \leq 1) \in \mathcal{C}[0,1]$$ is continuous $\mathbb{P}^{\bf W}$ a.s. According to Donsker's theorem \cite{Donsker}, see e.g. Billingsley \cite[Section $10$]{Bill} or Kallenberg \cite[Chapter $16$]{Kallenberg}, the linearly interpolated simple random walks
$$\left(\frac{RW_{[nt]}}{\sqrt{n}}; t \geq 0\right)~\mbox{converges weakly in $\mathcal{C}[0,1]$ to}~(B_t; t \geq 0),$$
So by the continuous mapping theorem, see e.g. Billingsley \cite[Theorem $5.1$]{Bill},
$$\Xi_F \circ \left(\frac{RW_{[nt]}}{\sqrt{n}}; t \geq 0\right)~\mbox{converges weakly to}~\Xi_F \circ (B_t; t \geq 0).~~\square$$
\quad Note that $T(\mathcal{BR}^{0,n})=F_n+n$. Following the above analysis, we know that $T(\mathcal{BR}^{0,n})/n$ converges weakly to $F+1$, where $T(\mathcal{BR}^{0,n})$ is the waiting time until an element of $\mathcal{BR}^{0,n}$ occurs in a simple random walk and $F$ is the random time defined as in \eqref{bridgelike}. As a consequence,
$$\liminf_{n \rightarrow \infty}\mathbb{E}\frac{T(\mathcal{BR}^{0,n})}{n} \geq (\mathbb{E}F+1),\quad \mbox{since}~\mathbb{E}F<\infty.$$
In particular, $\mathbb{E}F \leq C_{\mathcal{BR}}^0-1=3$ as in Section \ref{23}. We refer readers to Pitman and Tang \cite{PTacc} for further discussion on first level bridges and the structure of the Slepian zero set. 
\section{Continuous paths in Brownian motion}
\label{3}
\quad This section is devoted to the proof of Theorem \ref{thm3} and Theorem \ref{thm2}  regarding continuous paths and the distribution of continuous-time processes embedded in Brownian motion. In Section \ref{31}, we show that there is no normalized excursion in a Brownian path, Theorem \ref{thm2} \eqref{thm151}. A slight modification of the proof allows us to exclude the existence of the Vervaat bridges with negative endpoint, Theorem \ref{thm2} \eqref{thm155}. Furthermore, we prove in Section \ref{32} that there is even no reflected bridge in Brownian motion, Theorem \ref{thm2} \eqref{thm152}. In Section \ref{34} and \ref{35}, we show that neither the Vervaat transform of Brownian motion nor first passage bridges with negative endpoint can be found in Brownian motion, i.e. Theorem \ref{thm2} \eqref{thm153} \eqref{thm154}. We make use of the potential theory of {\em additive L\'{e}vy processes}, which is recalled in Section \ref{34}. Finally in Section \ref{33}, we provide a proof for the existence of Brownian meander, co-meander and three-dimensional Bessel process in Brownian motion, i.e. Theorem \ref{thm3}, using the acceptance-rejection method.
\subsection{No normalized excursion in a Brownian path}
\label{31}
In this part, we provide two proofs for \eqref{thm151} of Theorem \ref{thm2}, though similar, from different viewpoints. The first proof is based on a fluctuation version of {\em Williams' path decomposition} of Brownian motion, originally due to Williams \cite{Williams}, and later extended in various ways by Millar \cite{Millarbis,Millar}, and Greenwood and Pitman \cite{GreenPit}. We also refer readers to Pitman and Winkel \cite{PW} for a combinatorial explanation and various applications.
\begin{theorem} \cite{Williams,GreenPit}
\label{PW}
 Let $(B_t; t \geq 0)$ be standard Brownian motion and $\xi$ be exponentially distributed with rate $\frac{1}{2} \vartheta^2$, independent of $(B_t; t \geq 0)$. Define $M:=\argmin_{[0,\xi]}B_t$,  $H:=-B_M$ and $R:=B_{\xi}+H$. Then $H$ and $R$ are independent exponential variables, each with the same rate $\vartheta$. Furthermore, conditionally given $H$ and $R$, the path $(B_t; 0 \leq t \leq \xi)$ is decomposed into two independent pieces:
\begin{itemize}
\item
$(B_t; 0 \leq t \leq M)$ is Brownian motion with drift $-\vartheta<0$ running until it first hits the level $-H<0$;
\item
$(B_{\xi-t}-B_{\xi}; 0 \leq t \leq \xi-M)$ is Brownian motion with drift $-\vartheta<0$ running until it first hits the level $-R<0$.
\end{itemize}
\end{theorem}
\quad Now we introduce the notion of first passage process, which will be used in the proof of \eqref{thm151} of Theorem \ref{thm2}. Given a real-valued process $(Z_t; t \geq 0)$ starting at $0$, we define the first passage process $(\tau_{-x}; x \geq 0)$ associated to $X$ to be the first time that the level $-x < 0$ is hit: 
$$\tau_{-x}:=\inf \{t \geq 0; Z_t<-x\} \quad \mbox{for}~x > 0.$$
When $Z$ is Brownian motion, the distribution of the first passage process is well-known. It can be obtained either by random walks approximation or by {\em It\^{o}'s excursion theory} \cite{Itoex}.
\begin{lemma} 
\label{FPP} 
~
\begin{enumerate}
\item
\label{absl1}
Let ${\bf W}$ be Wiener measure on $\mathcal{C}[0,\infty)$. Then the first passage process $(\tau_{-x}; x \geq 0)$ under ${\bf W}$ is a stable$(\frac{1}{2})$ subordinator, with
$$\mathbb{E}^{{\bf W}}[\exp(-\alpha \tau_{-x})]=\exp(-x\sqrt{2 \alpha}) \quad \mbox{for}~\alpha>0.$$
\item
\label{absl2}

For $\vartheta \in \mathbb{R}$, let ${\bf W}^{\vartheta}$ be the distribution of Brownian motion with drift $\vartheta$. 

Then for each fixed $L>0$, on the event $\tau_{-L}<\infty$, the distribution of the first passage process $(\tau_{-x}; 0 \leq x \leq L)$ under ${\bf W}^{\vartheta}$ is absolutely continuous with respect to that under ${\bf W}$, with density $D^{\vartheta}_{L}:=\exp(-\vartheta L -\frac{\vartheta^2}{2} \tau_{-L})$.
\end{enumerate}
\end{lemma}
\textbf{Proof:} The part \eqref{absl1} of the lemma is a well known
result of L\'{e}vy. Or see Bertoin et al \cite[Lemma $4$]{BCP}. The part \eqref{absl2} is a direct consequence of Girsanov's theorem, see e.g. Revuz and Yor \cite[Chapter VIII]{RY} for background.  $\square$\\\\
\textbf{Proof of Theorem \ref{thm2} \eqref{thm151}:} Suppose by contradiction that $\mathbb{P}(T < \infty)>0$, where $T$ is the first time that some excursion pattern appears. Take $\xi$ exponentially distributed with rate $\frac{1}{2}$, independent of $(B_t; t \geq 0)$. We have then 
\begin{equation}
\label{1}
\mathbb{P}(T < \xi <T+1)>0.
\end{equation}
Now $(T,T+1)$ is inside the excursion of Brownian motion above its past-minimum process, which straddles $\xi$. Define
\begin{itemize}
\item
$(\sigma_{-x}; x \geq 0)$ to be the first passage process of $(B_{\xi-t}-B_{\xi}; 0 \leq t \leq \xi-M)$ concatenated by an independent Brownian motion with drift $-1$ running forever ;
\item
$(\tau_{-x}; x \geq 0)$ to be the first passage process of $(B_{\xi+t}-B_{\xi}; t \geq 0)$.
\end{itemize}
\begin{figure}[h]
\includegraphics[width=0.6 \textwidth]{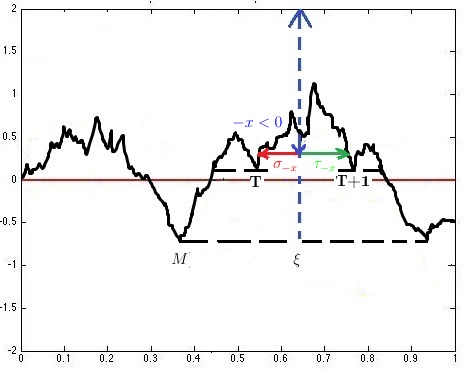}
\caption{No excursion of length $1$ in a Brownian path.}
\end{figure}
By the strong Markov property of Brownian motion, $(B_{\xi+t}-B_{\xi}; t \geq 0)$ is still Brownian motion. Thus, $(\tau_{-x}; x \geq 0)$ is a stable$(\frac{1}{2})$ subordinator by part \eqref{absl1} of Lemma \ref{FPP}. According to Theorem \ref{PW}, $(B_{\xi-t}-B_{\xi}; 0 \leq t \leq \xi-M)$ is Brownian motion with drift $-1$ running until it first hits the level $-R<0$. Then $(\sigma_{-x}; x \geq 0)$ is the first passage process of Brownian motion with drift $-1$, whose distribution is absolutely continuous on any compact interval $[0,L]$, with respect to that of $(\tau_{-x}; 0 \leq x \leq L)$ by part \eqref{absl2} of Lemma \ref{FPP}. It is well known that a real stable$(\frac{1}{2})$ process does not hit points, see e.g. Bertoin \cite[Theorem $16$, Chapter II.$5$]{Bertoin}. As a consequence,
$$\mathbb{P}(\sigma_{-x}+\tau_{-x}=1~\mbox{for some}~x \geq 0)=0,$$
which contradicts \eqref{1}.  $\square$

\quad It is easy to see that the above argument still works for proving the non-existence of the Vervaat bridge pattern $\mathcal{VB}^{\lambda}$, with endpoint $\lambda<0$.\\\\
\textbf{Proof of Theorem \ref{thm2} \eqref{thm155}:} We borrow the notations from the preceding proof. Observe that, for fixed $\lambda<0$,
$$\mathbb{P}(\sigma_{-x}+\tau_{-x+\lambda}=1~\mbox{for some}~x \geq 0)=0.$$ The rest of the proof is just a duplication of the preceding one.  $\square$

\quad We give yet another proof of Theorem \ref{thm2} \eqref{thm151}, which relies on It\^{o}'s excursion theory, combined with Bertoin's self-similar fragmentation theory. For general background on fragmentation processes, we refer to the monograph of Bertoin \cite{BertoinFC}. The next result, regarding a normalized Brownian excursion, follows Bertoin \cite[Corollary $2$]{selfBertoin}.
\begin{theorem} \cite{selfBertoin}
\label{selfsimilar}
Let $e:=(e_u; 0 \leq u \leq 1)$ be normalized Brownian excursion and $F^{e}:=(F_t^{e}; t \geq 0)$ be the associated interval fragmentation defined as $F^{e}_t:=\{u \in (0,1); e_u>t\}$. Introduce
\begin{itemize}
\item
$\lambda:=(\lambda_t; t \geq 0)$ the length of the interval component of $F^{e}$ that contains $U$, independent of the excursion and uniformly distributed;
\item
$\xi:=\{\xi_t; t \geq 0\}$ a subordinator, the Laplace exponent of which is given by
\begin{equation}
\label{subord}
\Phi^{ex}(q):=q \sqrt{\frac{8}{\pi}} \int_{0}^1t^{q-\frac{1}{2}}(1-t)^{-\frac{1}{2}}= q \sqrt{\frac{8}{\pi}} B(q+\frac{1}{2},\frac{1}{2});
\end{equation}
\end{itemize}
Then $(\lambda_t;t \geq 0)$ has the same law as $(\exp(-\xi_{\rho_t}); t \geq 0)$, where 
\begin{equation}
\label{timechange}
\rho_t:=\inf \left\{u \geq 0; \int_0^u \exp \left(-\frac{1}{2}\xi_r\right)dr>t \right\}.
\end{equation}
\end{theorem}
\textbf{Alternative proof of Theorem \ref{thm2} \eqref{thm151}:} Consider the reflected process $(B_t-\underline{B}_t; t \geq 0)$, where $\underline{B}_t:=\inf_{0 \leq u \leq t}B_u$ is the past-minimum process of the Brownian motion. For $\mathbf{e}$ the first excursion of $B-\underline{B}$ that contains some excursion pattern $\mathcal{E}$ of length $1$, let $\Lambda_{\mathbf{e}}$ be the length of such excursion, and $\mathbf{e}^{*}$ be the normalized Brownian excursion. Following It\^{o}'s excursion theory, see e.g. Revuz and Yor \cite[Chapter XII]{RY}, $\Lambda_{\mathbf{e}}$ is independent of the distribution of the normalized excursion $\mathbf{e}^{*}$.  As a consequence, the fragmentation associated to $\mathbf{e}^{*}$ produces an interval of length $\frac{1}{\Lambda_{\mathbf{e}}}$. Now choose $U$ uniformly distributed on $[0,1]$ and independent of the Brownian motion.  According to Theorem \ref{selfsimilar}, there exists a subordinator $\xi$ characterized as in \eqref{subord} and a time-change $\rho$ defined as in \eqref{timechange} such that $(\lambda_t; t \geq 0)$, the process of the length of the interval fragmentation which contains $U$, has the same distribution as $(\exp(-\xi_{\rho_t});t \geq 0)$.
Note that $(\lambda_t;t \geq 0)$ depends only on the normalized excursion $\mathbf{e}^{*}$ and $U$, so $(\lambda_t;t \geq 0)$
is independent of $\Lambda_{\mathbf{e}}$. It is a well known result of Kesten \cite{Kesten69} that a subordinator without drift does not hit points. Therefore,
$$\mathbb{P}\left(\lambda_t=\frac{1}{\Lambda_{\mathbf{e}}}~\mbox{for some}~t \geq 0\right)=0,$$
which yields the desired result.  $\square$
\subsection{No reflected bridge in a Brownian path}
\label{32}
This part is devoted to proving Theorem \ref{thm2} \eqref{thm152}. The main difference between Theorem \ref{thm2} \eqref{thm151} and \eqref{thm152} is that the strict inequality $B_{T+u}>B_{T}$ for all $u \in (0,1)$ is relaxed by the permission of equalities $B_{T+u}=B_T$ for some $u \in (0,1)$. Thus, there seems to be more chance to find the reflected bridge paths $\mathcal{RBR}$ in Brownian motion. Nevertheless, the following lemma suggests that one cannot expect the equality to be achieved for more than one instant in between. Below is a slightly stronger version of this result.
\begin{lemma}
\label{twice}
Almost surely, there are no random times $S<T$ such that $B_{T}= B_{S}$, $B_{u} \geq B_{S}$ for $u \in (S,T)$ and $B_{v}=B_{w}=B_S$ for some $S<v<w<T$.
\end{lemma}
\textbf{Proof:} Consider the following two sets $$\mathcal{T}:=\{\mbox{there exist}~S~\mbox{and}~T~\mbox{which satisfy the conditions in the lemma}\}$$ and $$\mathcal{U}:=\bigcup_{s, t \in \mathbb{Q}}\{B~\mbox{attains its minimum for more than once on}~[s,t]\}.$$
\quad It is straightforward that $\mathcal{T} \subset \mathcal{U}$. In addition, it is well-known that almost surely Brownian motion has a unique minimum on any fixed interval $[s,t]$ for all $s,t \in \mathbb{R}$. As a countable union of null events, $\mathbb{P}(\mathcal{U})=0$ and thus $\mathbb{P}(\mathcal{T})=0$.  $\square$
\begin{remark}
{\em The previous lemma has an interesting geometric interpretation in terms of Brownian trees, see e.g. Pitman \cite[Section $7.4$]{Pitman} for background. Along the lines of the second proof of Theorem \ref{thm2} \eqref{thm151} in Section \ref{31}, we only need to show that the situation in Lemma \ref{twice} cannot happen in Brownian excursion either of an independent and diffuse length or of normalized unit length. But this is just another way to state that Brownian trees have only binary branch points, which follows readily from Aldous' stick-breaking construction of the continuum random trees, see e.g. Aldous \cite[Section $4.3$]{Aldous3} and Le Gall \cite{LeGall}.}
\end{remark}
\quad According to Theorem \ref{thm2} \eqref{thm151} and Lemma \ref{twice}, we see that almost surely, there are neither excursion paths nor reflected bridge paths with at least two reflections in Brownian motion. To prove the desired result, it suffices to exempt the possibility of reflected bridge paths with exactly one reflection.
\begin{lemma}
\label{once}
Assume that $0 \leq S<T<U$ are random times such that $B_S=B_T=B_U$ and $B_u>B_S$ for $u \in (S,T) \cup (T,U)$. Then the distribution of $U-S$ is absolutely continuous with respect to the Lebesgue measure.
\end{lemma}
\textbf{Proof:} Suppose by contradiction that the distribution of $U-S$ is not absolutely continuous with respect to the Lebesgue measure. Then there exists $p,q \in \mathbb{Q}$ such that $U-S$ fails to have a density on the event $\{S<p<T<q<U\}$. In fact,  if $U-S$ has a density on $\{S<p<T<q<U\}$ for all $p,q \in \mathbb{Q}$, Radon-Nikodym theorem guarantees that $U-S$ has a density on $\{S<T<U\}=\cup_{p,q \in \mathbb{Q}}\{S<p<T<q<U\}$. Note that on the event $\{S<p<T<q<U\}$, $U$ is the first time after $q$ such that the Brownian motion $B$ attains $\inf_{u \in [p,q]} B_u$ and obviously has a density. Again by Radon-Nikodym theorem, the distribution of $U-S$ has a density on $\{S<p<T<q<U\}$, which leads to a contradiction.  $\square$
\begin{remark}
{\em The previous result can also be inferred from a fine study on local minima of Brownian motion. Neveu and Pitman \cite{NeveuPitman} studied the renewal structure of local extrema in a Brownian path, in terms of Palm measure, see e.g. Kallenberg \cite[Chapter $11$]{Kallenberg} for background. More precisely, denote 
\begin{itemize}
\item 
$\mathcal{C}$ to be the space of continuous paths on $\mathbb{R}$, equipped with Wiener measure ${\bf W}$;
\item
$E$ to be the space of excursions with lifetime $\zeta$, equipped with It\^{o} measure ${\bf n}$.
\end{itemize}
Then the Palm measure of all local minima is the image of $\frac{1}{2}({\bf n} \times {\bf n} \times {\bf W})$ by the mapping $E \times E \times \mathcal{C} \ni (e,e',w) \rightarrow  \tilde{w} \in \mathcal{C}$ given by 
$$\tilde{w}_t =     \left\{ \begin{array}{ccl}
         w_{t+\zeta(e')} & \mbox{if}
         & t \leq -\zeta(e'), \\ e'_{-t}  & \mbox{if} & -\zeta(e') \leq t \leq 0, \\
         e_t & \mbox{if} & 0 \leq t \leq \zeta(e), \\ w_{t-\zeta(e)} & \mbox{if} & t \geq \zeta(e).
                \end{array}\right.$$
\begin{figure}[h]
\includegraphics[width=0.6 \textwidth]{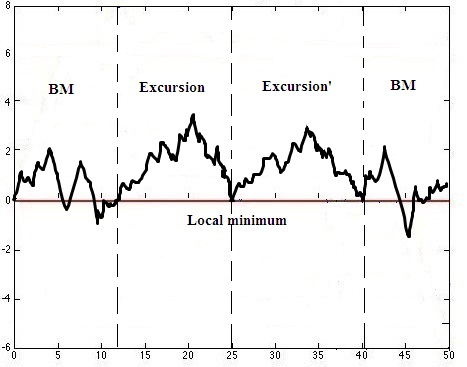}
\caption{Structure of local minima in Brownian motion.}
\end{figure}

Using the notations of Lemma \ref{once}, an in-between reflected position $T$ corresponds to a Brownian local minimum. Then the above discussion implies that $U-S$ is the sum of two independent random variables with densities and hence is diffuse. See also Tsirelson \cite{Tsirelson} for the i.i.d. uniform sampling construction of Brownian local minima, which reveals the diffuse nature of $U-S$.}
\end{remark}
\subsection{No Vervaat tranform of Brownian motion in a Brownian path}
\label{34}
In the current section, we aim to prove Theorem \ref{thm2} \eqref{thm153}. That is, there is no random time $T$ such that $$(B_{T+u}-B_T; 0 \leq u \leq 1) \in \mathcal{VB}^{-},$$
where $\mathcal{VB}^{-}:=\{w \in \mathcal{C}[0,1]; w(t)>w(1)~\mbox{for}~0 \leq t<1~\mbox{and}~\inf\{t>0;w(t)<0\}>0\}$. A similar argument shows that there is no random time $T$ such that $$(B_{T+u}-B_T; 0 \leq u \leq 1) \in \mathcal{VB}^{+},$$
where $\mathcal{VB}^{+}:=\{w \in \mathcal{C}[0,1]; w(t)>w(0)~\mbox{for}~0 < t \leq 1~\mbox{and}~\sup\{t<1;w(t)<w(1)\}>0\}$. Observe that $\mathcal{VB}^{-} \cup \mathcal{VB}^{-}$ is the support of $(V_u; 0 \leq u \leq 1)$. Thus, the Vervaat transform of Brownian motion cannot be embedded into Brownian motion.

\quad In Section \ref{31}, we showed that for each fixed $\lambda<0$, there is no random time $T$ such that $(B_{T+u}-B_{T}; 0 \leq u \leq 1) \in \mathcal{VB}^{\lambda}$. However, there is no obvious way to pass from the non-existence of the Vervaat bridges to that of the Vervaat transform of Brownian motion, due to an uncountable number of possible levels. 

\quad To get around the problem, we make use of an additional tool -- potential theory of additive L\'{e}vy processes, developed by Khoshnevisan et al \cite{Khosh2002,Khosh2003c,Khosh2003b,Khosh2003a,Khosh2009}. We now recall some results of this theory that we need in the proof of Theorem \ref{thm2} \eqref{thm153}. For a more extensive overview of the theory, we refer readers to the survey of Khoshnevisan and Xiao \cite{Khosh2005}.
\begin{definition}
\label{adlp}
An $N$-parameter, $\mathbb{R}^d$-valued additive L\'{e}vy process $(Z_{\textbf{t}}; \textbf{t} \in \mathbb{R}_{+}^N)$ with L\'{e}vy exponent $(\Psi^1,\ldots,\Psi^N)$ is defined as 
\begin{equation}
\label{alp}
Z_{\textbf{t}}:=\sum_{i=1}^N Z^i_{t_i} \quad \mbox{for}~\textbf{t}=(t_1,\ldots,t_N) \in \mathbb{R}_{+}^N,
\end{equation}
where $(Z^1_{t_1};t_1 \geq 0),\ldots,(Z^N_{t_N};t_N \geq 0)$ are $N$ independent $\mathbb{R}^d$-valued L\'{e}vy processes with L\'{e}vy exponent $\Psi^1,\ldots,\Psi^N$.
\end{definition}
\quad The following result regarding the range of additive L\'{e}vy processes is due to Khoshnevisan et al \cite[Theorem $1.5$]{Khosh2003a}, \cite[Theorem $1.1$]{Khosh2009}, and Yang \cite[Theorem $1.1$]{Yang2007,Yang2009}.
\begin{theorem} \cite{Khosh2003a,Yang2007,Khosh2009}
\label{KXY}
Let $({Z_{\bf t}; {\bf t} \in \mathbb{R}_{+}^N})$ be an additive L\'{e}vy process defined as in \eqref{alp}. Then
$$\mathbb{E}[\leb(Z(\mathbb{R}_{+}^N))]>0 \Longleftrightarrow \int_{\mathbb{R}^d} \prod_{i=1}^N \re\left(\frac{1}{1+\Psi^i(\zeta)}\right)d\zeta<\infty,$$
where $\leb(\cdot)$ is the Lebesgue measure on $\mathbb{R}^d$, and $\re(\cdot)$ is the real part of a complexe number.
\end{theorem}
\quad The next result, which is read from Khoshnevisan and Xiao \cite[Lemma $4.1$]{Khosh2005b}, makes a connection between the range of an additive L\'{e}vy process and the polarity of single points. See also Khoshnevisan and Xiao \cite[Lemma 3.1]{Khosh2005}.
\begin{theorem} \cite{Khosh2002,Khosh2005b}
\label{KX}
Let $({Z_{\bf t}; {\bf t} \in \mathbb{R}_{+}^N})$ be an additive L\'{e}vy process defined as in \eqref{alp}. Assume that for each $\textbf{t} \in \mathbb{R}_{+}^{N}$, the distribution of $Z_{\textbf{t}}$ is mutually absolutely continuous with respect to Lebesgue measure on $\mathbb{R}^d$. Let $z \in \mathbb{R}^d \setminus \{0\}$, then
$$\mathbb{P}(Z_{\textbf{t}}=z~\mbox{for some}~\textbf{t} \in \mathbb{R}_{+}^N)>0 \Longleftrightarrow \mathbb{P}(\leb(Z(\mathbb{R}_{+}^N)>0)>0.$$
\end{theorem}
\quad Note that $\mathbb{P}(\leb(Z(\mathbb{R}_{+}^N)>0)>0$ is equivalent to $\mathbb{E}[\leb(Z(\mathbb{R}_{+}^N))]>0$. Combining Theorem \ref{KXY} and Theorem \ref{KX}, we have:
\begin{corollary}
\label{imp}
Let $({Z_{\bf t}; {\bf t} \in \mathbb{R}_{+}^N})$ be an additive L\'{e}vy process defined as in \eqref{alp}. Assume that for each $\textbf{t} \in \mathbb{R}_{+}^{N}$, the distribution of $Z_{\textbf{t}}$ is mutually absolutely continuous with respect to Lebesgue measure on $\mathbb{R}^d$.  Let $z \in \mathbb{R}^d \setminus \{0\}$, then
$$\mathbb{P}(Z_{\textbf{t}}=z~\mbox{for some}~\textbf{t} \in \mathbb{R}_{+}^N)>0 \Longleftrightarrow \int_{\mathbb{R}^d} \prod_{i=1}^N \re\left(\frac{1}{1+\Psi^i(\zeta)}\right)d\zeta<\infty.$$
\end{corollary}
\textbf{Proof of Theorem \ref{thm2} \eqref{thm153}:} We borrow the notations from the proof of Theorem \ref{thm2} \eqref{thm151} in Section \ref{31}. It suffices to show that 
\begin{equation}
\label{bber}
\mathbb{P}(\sigma_{-t_1}+\tau_{-t_2}=1 ~\mbox{for some}~t_1,t_2 \geq 0)=0,
\end{equation}
where $(\sigma_{-t_1};t_1 \geq 0)$ is the first passage process of Brownian motion with drift $-1$, and $(\tau_{-t_2};t_2 \geq 0)$ is a stable$(\frac{1}{2})$ subordinator independent of $(\sigma_{-t_1};t_1 \geq 0)$. Let $Z_{\bf t}=Z_{t_1}^1+Z_{t_2}^2:=\sigma_{-t_1}+\tau_{-t_2}$ for ${\bf t}=(t_1,t_2) \in \mathbb{R}_{+}^2$. By Definition \ref{adlp}, $Z$ is a $2$-parameter, real-valued additive L\'{e}vy process with L\'{e}vy exponent $(\Psi^1,\Psi^2)$ given by
$$\Psi^1(\zeta)=\sqrt[4]{1+4\zeta^2}\exp\left[-i\frac{\arctan(2\zeta)}{2}\right]-1 \quad \mbox{and} \quad \Psi^2(\zeta)=\sqrt{|\zeta|}(1-i\sgn \zeta) \quad \mbox{for}~\zeta \in \mathbb{R},$$
which is derived from the formula in Cinlar \cite[Chapter $7$, Page $330$]{Cinlar} and Lemma \ref{FPP} \eqref{absl2}. Hence,
$$\re\left(\frac{1}{1+\Psi^1(\zeta)}\right)=\frac{1}{\sqrt[4]{1+4\zeta^2}}\sqrt{\frac{1}{2}\left(1+\frac{1}{\sqrt{1+4\zeta^2}}\right)} ~\mbox{and} ~ \re\left(\frac{1}{1+\Psi^2(\zeta)}\right)=\frac{1+\sqrt{|\zeta}|}{1+2\sqrt{|\zeta|}+2|\zeta|}.$$
Clearly, $\Xi: \zeta \rightarrow \re\left(\frac{1}{1+\Psi^1(\zeta)}\right)\re\left(\frac{1}{1+\Psi^2(\zeta)}\right)$ is not integrable on $\mathbb{R}$ since $\Xi(\zeta) \sim \frac{1}{4 |\zeta|}$ as $|\zeta| \rightarrow \infty$. In addition, for each ${\bf t} \in \mathbb{R}_{+}^2$, $Z_t$ is absolutely continuous with respect to Lebesgue measure on $\mathbb{R}$. Applying Corollary \ref{imp}, we obtain \eqref{bber}.  $\square$
\subsection{No first passage bridge in a Brownian path}
\label{35}
We prove Theorem \ref{thm2} \eqref{thm154}, i.e. there is no first passage bridge in Brownian motion by a spacetime shift. The main difference between Vervaat bridges with fixed endpoint $\lambda<0$ and first passage bridges ending at $\lambda<0$ is that the former start with an excursion piece, while the latter return to the origin infinitely often on any small interval $[0,\epsilon]$, $\epsilon>0$. Thus, the argument used in Section \ref{31} to prove the non-existence of Vervaat bridges is not immediately applied in case of first passage bridges. Nevetheless, the potential theory of additive L\'{e}vy processes helps to circumvent the difficulty.\\\\
\textbf{Proof of Theorem \ref{thm2} \eqref{thm154}:} Suppose by contradiction that $\mathbb{P}(T < \infty)>0$, where $T$ is the first time that some first passage bridge appears. Take $\xi$ exponentially distributed with rate $\frac{1}{2}$, independent of $(B_t; t \geq 0)$. We have then 
\begin{equation}
\label{1}
\mathbb{P}(T < \xi <T+1)>0.
\end{equation}
Now $(T,T+1)$ is inside the excursion of Brownian motion below its past-maximum process, which straddles $\xi$. Define $(\tau_{-x}; x \geq 0)$ to be the first passage process of $(B_{\xi+t}-B_{\xi}; t \geq 0)$.
\begin{figure}[h]
\includegraphics[width=0.6 \textwidth]{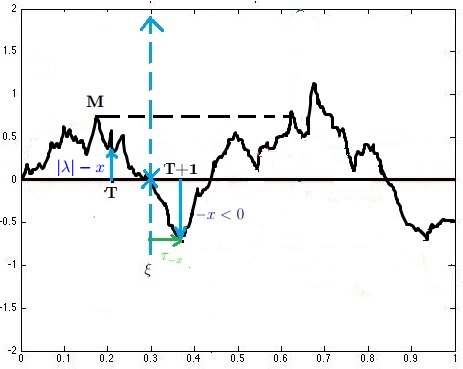}
\caption{No first passage bridge of length $1$ in a Brownian path.}
\end{figure}

By strong Markov property of Brownian motion, $(B_{\xi+t}-B_{\xi}; t \geq 0)$ is still Brownian motion. Thus, $(\tau_{-x}; x \geq 0)$ is a stable$(\frac{1}{2})$ subordinator. Let $M:=\argmax_{[0,\xi]}B_t$. By a variant of Theorem \ref{PW}, $(B_{\xi-t}-B_{\xi}; 0 \leq t \leq \xi-M)$ is Brownian motion with drift $1$ running until it first hits the level $B_M-B_{\xi}>0$, independent of $(\tau_{-x}; x \geq 0)$. As a consequence,  \eqref{1} implies that 
\begin{equation}
\label{11}
\mathbb{P}(\tau_{-x}=l ~\mbox{and}~B^{\uparrow}_{1-l}=|\lambda|-x~\mbox{for some}~(x,l) \in \mathbb{R}_{+} \times [0,1])>0,
\end{equation}
where $(B^{\uparrow}_t; t \geq 0)$ is Brownian motion with drift $1$, independent of $\frac{1}{2}$-stable subordinator $(\tau_{-x}; x \geq 0)$. By setting $t_1:=x$ and $t_2:=1-l$, we have:
\begin{align}
&~\quad \mathbb{P}(\tau_{-x}=l ~\mbox{and}~B^{\uparrow}_{1-l}=|\lambda|-x~\mbox{for some}~(x,l) \in \mathbb{R}_{+} \times [0,1])\notag\\
&=\mathbb{P}(\tau_{-t_1}+t_2=1 ~\mbox{and}~B^{\uparrow}_{t_2}+t_1=|\lambda|~\mbox{for some}~(t_1,t_2) \in \mathbb{R}_{+} \times [0,1]) \notag\\
& \leq \mathbb{P}[(\tau_{-t_1},t_1)+(t_2,B^{\uparrow}_{t_2})=(1,|\lambda|)~\mbox{for some}~(t_1,t_2) \in \mathbb{R}_{+}^2] \label{22}
\end{align}
Let $Z_{\bf t}=Z^1_{t_1}+Z^2_{t_2}:=(\tau_{-t_1},t_1)+(t_2,B^{\uparrow}_{t_2})$ for $\textbf{t}=(t_1,t_2) \in \mathbb{R}_{+}^2$. By Definition \ref{adlp}, $Z$ is a $2$-parameter, $\mathbb{R}^2$-valued additive L\'{e}vy process with L\'{e}vy exponent $(\Psi^1,\Psi^2)$ given by
$$\Psi^1(\zeta_1,\zeta_2):=\sqrt{|\zeta_1|}-i(\sqrt{|\zeta_1|} \sgn \zeta_1+\zeta_2) \quad \mbox{and} \quad \Psi^2(\zeta_1,\zeta_2):=\frac{\zeta_2^2}{2}-i(\zeta_1+\zeta_2) \quad \mbox{for}~(\zeta_1,\zeta_2) \in \mathbb{R}^2.$$
Hence,
\begin{multline*}
\re\left(\frac{1}{1+\Psi^1(\zeta_1,\zeta_2)}\right)\re\left(\frac{1}{1+\Psi^2(\zeta_1,\zeta_2)}\right)\\
=\frac{(1+\sqrt{|\zeta_1|})\left(1+\frac{\zeta_2^2}{2}\right)}{\left[(1+\sqrt{|\zeta_1|})^2+(\sqrt{|\zeta_1|} \sgn \zeta_1+\zeta_2)^2 \right] \left[\left(1+\frac{\zeta_2^2}{2}\right)^2+(\zeta_1+\zeta_2)^2\right]}:=\Xi(\zeta_1,\zeta_2). 
\end{multline*}
Observe that $\zeta \rightarrow \Xi(\zeta_1,\zeta_2)$ is not integrable on $\mathbb{R}^2$, which is clear by passage to polar coordinates $(\zeta_1,\zeta_2)=(\rho \cos \theta, \sqrt{\rho}\sin \theta)$ for $\rho \geq 0$, $\theta \in [0,2 \pi)$. In addition, for each $\textbf{t} \in \mathbb{R}_{+}^2$, $Z_{\bf t}$ is mutually absolutely continuous with respect to Lebesgue measure on $\mathbb{R}^2$. Applying Corollary \ref{imp}, we know that
$$\mathbb{P}(Z_{\textbf{t}}=(1,|\lambda|)~\mbox{for some}~\textbf{t} \in \mathbb{R}_{+}^2)=0.$$
Combining with \eqref{22}, we obtain:
$$\mathbb{P}(\tau_{-x}=l ~\mbox{and}~B^{\uparrow}_{1-l}=|\lambda|-x~\mbox{for some}~(x,l) \in \mathbb{R}_{+} \times [0,1])=0,$$
which contradicts \eqref{11}.  $\square$

\quad It is well-known that time reversal of Brownian first passage bridge ending at $\lambda<0$, i.e.
$$(F^{\lambda,br}_{1-u}+|\lambda|; 0 \leq u\leq 1)$$
has the same distribution as three dimensional Bessel bridge ending at $|\lambda|>0$, see e.g. Biane and Yor \cite{BY}. It is not hard to see that the above argument still works for $\delta$-dimensional Bessel bridges with $\delta \geq 2$, where $\{0\}$ is polar.
\begin{corollary} (Non-existence of positive bridge paths/$\delta$-dimensional Bessel bridge, $\delta \geq 2$) For each fixed $\lambda>0$, almost surely, there is no random time $T$ such that 
$$(B_{T+u}-B_T; 0 \leq u \leq 1) \in \mathcal{BES}^{\lambda}:=\{w \in \mathcal{C}[0,1]; w(t)>0~\mbox{and}~w(1)=\lambda\}.$$
In particular, there is no random time $T \geq 0$ such that $(B_{T+u}-B_T; 0 \leq u \leq 1)$ has the same distribution as $\delta$-dimensional Bessel bridge ending at $\lambda$, with $\delta \geq 2$.
\end{corollary}
\subsection{Meander, co-meander and 3-d Bessel process in a Brownian path}
\label{33}
We prove Theorem \ref{thm3} in this section using It\^{o}'s excursion theory, combined with {\em Rost's filling scheme} \cite{CO,Rost1} solution to the Skorokhod embedding problem.

\quad The existence of Brownian meander in a Brownian path is assured by the following well-known result, which can be read from Maisoneuve \cite[Section $8$]{Maison}, with explicit formulas due to Chung \cite{Chung}. An alternative approach was provided by Greenwood and Pitman \cite{GreenPitbis}, and Pitman \cite[Section $4$ and $5$]{Pitmanlevy}. See also Biane and Yor \cite[Theorem $6.1$]{BianeYor}, or Revuz and Yor \cite[Exercise $4.18$, Chapter XII]{RY}.
\begin{theorem} \cite{Maison,Pitmanlevy, BianeYor}
\label{meander}
Let $(e^i)_{i \in \mathbb{N}}$ be the sequence of excursions, whose length exceeds $1$, in the reflected process $(B_t-\underline{B}_t; t \geq 0)$. Then $(e^i_u;0 \leq u \leq 1)_{i \in \mathbb{N}}$ is a sequence of independent and identically distributed paths, each distributed as Brownian meander $(m_u; 0 \leq u \leq 1)$.
\end{theorem}
\quad Let us recall another basic result due to Imhof \cite{Imhof}, which establishes the absolute continuity relation between Brownian meander and the three-dimensional Bessel process. Their relation with Brownian co-meander is studied in Yen and Yor \cite[Chapter $7$]{YY}.
\begin{theorem} \cite{Imhof, YY} \label{abscond} The distributions of Brownian meander $(m_u;0 \leq u \leq 1)$, Brownian co-meander $(\widetilde{m}_u; 0 \leq u \leq 1)$ and the three-dimensional Bessel process $(R_u; 0 \leq u \leq 1)$ are mutually absolutely continuous with respect to each other. For $F: \mathcal{C}[0,1] \rightarrow \mathbb{R}^{+}$ a measurable function, 
\begin{enumerate}
\item
$\mathbb{E}[F(m_u;0 \leq u \leq 1)]=\mathbb{E}\left[\sqrt{\frac{\pi}{2}} \frac{1}{R_1} F(R_u;0 \leq u \leq 1)\right]$;
\item
$\mathbb{E}[F(\widetilde{m}_u;0 \leq u \leq 1)]=\mathbb{E}\left[\frac{1}{R_1^2} F(R_u;0 \leq u \leq 1)\right]$.
\end{enumerate}
\end{theorem}
\quad According to Theorem \ref{meander}, there exist $T_1,T_2, \cdots$ such that 
\begin{equation}
\label{iidme}
m^i:=(B_{T_i+u}-B_{T_i}; 0 \leq u \leq 1)
\end{equation} 
form a sequence of i.i.d. Brownian meanders. Since Brownian co-meander and the three-dimensional Bessel process are absolutely continuous relative to Brownian meander, it is natural to think of {\em von Neumann's acceptance-rejection algorithm} \cite{vN}, see e.g. Rubinstein and Kroese \cite[Section $2.3.4$]{Rub} for background and various applications. However, von Neumann's selection method requires that the Radon-Nikodym density between the underlying probability measures is essentially bounded, which is not satisfied in the cases suggested by Theorem \ref{abscond}. Nevertheless, we can apply the filling scheme of Chacon and Ornstein \cite{CO} and Rost \cite{Rost1}.

\quad We observe that sampling Brownian co-meander or the three-dimensional Bessel process from i.i.d. Brownian meanders $(m^i)_{i\in \mathbb{N}}$ fits into the general theory of Rost's filling scheme applied to the Skorokhod embedding problem. In the sequel, we follow the approach of Dellacherie and Meyer \cite[Section $63-74$, Chapter IX]{DMN}, which is based on the seminal work of Rost \cite{Rost1}, to construct a stopping time $N$ such that $m^N$ achieves the distribution of $\widetilde{m}$ or $R$. We need some notions from potential theory for the proof.
\begin{definition} \label{defp}
~
\begin{enumerate}
\item 
\label{def1}
Given a Markov chain $X:=(X_n)_{n \in \mathbb{N}}$, a function $f$ is said to be excessive relative to $X$ if $$(f(X_n))_{n \in \mathbb{N}}~\mbox{is}~\mathcal{F}_n-\mbox{supermartingale},$$ where $(\mathcal{F}_n)_{n \in \mathbb{N}}$ is the natrual filtrations of $X$.
\item 
\label{def2}
Given two positive measures $\mu$ and $\lambda$, $\mu$ is said to be a balayage/sweeping of $\lambda$ if
$$\mu(f) \leq \lambda(f) \quad \mbox{for all bounded excessive functions}~f.$$
\end{enumerate}
\end{definition}
\textbf{Proof of Theorem \ref{thm3}:} Let $\mu^m$ (resp. $\mu^R$) be the distribution of Brownian meander (resp. the three-dimensional Bessel process) on the space $(\mathcal{C}[0,1],\mathcal{F})$. By the filling scheme, the sequence of measures $(\mu^m_i,\mu^R_i)_{i \in \mathbb{N}}$ is defined recursively as
\begin{equation}
\label{fs0}
\mu^m_0:=(\mu^m-\mu^R)^{+}\quad \mbox{and} \quad \mu^R_0:=(\mu^m-\mu^R)^{-},
\end{equation}
and for each $i \in \mathbb{N}$,
\begin{equation}
\label{fsi}
\mu^m_{i+1}:=(\mu^m_i(1)\cdot \mu^m-\mu^R_i)^{+}\quad \mbox{and} \quad \mu^R_{i+1}:=(\mu^m_i(1)\cdot \mu^m-\mu^R_i)^{-},
\end{equation}
where $\mu^m_i(1)$ is the total mass of the measure $\mu^m_i$. It is not hard to see that the bounded excessive functions of the i.i.d. meander sequence are constant $\mu^m$ a.s. Since $\mu^R$ is absolutely continuous with respect to $\mu^m$, for each $\mu^m$ a.s. constant function $c$, $\mu^R(c)=\mu^m(c)=c$.
Consequently,  $\mu^R$ is a balayage/sweeping of $\mu^m$ by Definition \ref{defp}. According to Theorem $69$ of Dellacherie and Meyer \cite{DMN}, $$\mu^R_{\infty}=0, \quad \mbox{where}~\mu^R_{\infty}:=\downarrow \lim_{i \rightarrow \infty}\mu^R_i.$$  

Now let $d_0$ be the Radon-Nikodym density of $\mu^m_0$ relative to $\mu^m$, and for $i>0$, $d_i$ be the Radon-Nikodym density of $\mu^m_i$ relative to $\mu^m_{i-1}(1)\cdot \mu^m$. We have
\begin{align}
\label{intem}
\mu^R &= (\mu^R-\mu^R_0)+(\mu^R_0-\mu^R_1)+\cdots \notag\\
  &= (\mu^m-\mu^m_0)+(\mu^m_0(1) \cdot \mu^m-\mu^m_1)+\cdots \notag\\
  &= (1-d_0)\mu^m+d_0\mu^m(1) \cdot (1-d_1)\mu^m+\cdots. 
\end{align}
Consider the stopping time $N$ defined by
\begin{equation}
\label{Nstop}
N:=\inf\left\{n \geq 0; -\sum_{i=0}^n \log d_i(m^i)>\xi\right\},
\end{equation}
where $(d_i)_{i \in \mathbb{N}}$ is the sequence of Radon-Nikodym densities defined as in the preceding paragraph, $(m^i)_{i \in \mathbb{N}}$ is the sequence of i.i.d. Brownian meanders defined as in \eqref{iidme}, and $\xi$ is exponentially distributed with rate $1$, independent of $(m^i)_{i \in \mathbb{N}}$. 

From the computation of \eqref{intem}, for all bounded measurable function $f$ and all $k \in \mathbb{N}$,
\begin{align*}
\mathbb{E}[f(m^N);N=k] &= \mathbb{E}[f(m^k);-\sum_{i=0}^{k-1} \log d_i(m^i) \leq \xi <-\sum_{i=0}^k \log d_i(m^i)] \\
                       &= \mathbb{E}[d_0(m^0) \cdots d_{k-1}(m^{k-1})f(m^k)(1-d_k(m^k))] \\
                       &= (\mu^m_{k-1}(1) \cdot \mu^m -\mu^m_k)f \\
                       &= (\mu^R_{k-1}-\mu^R_k)f,
\end{align*}
where $(\mu^m_i, \mu^R_i)_{i \in \mathbb{N}}$ are the filling measures defined as in \eqref{fs0} and \eqref{fsi}. By summing over all $k$, we have
$$\mathbb{E}[f(m^N);N<\infty]=\mu^R f,$$
i.e. $m^N$ has the same distribution as $R$. As a summary, 
$$(B_{T_N+u}-B_{T_N}; 0 \leq u \leq 1)~\mbox{has the same distribution as}~(R_u; 0 \leq u \leq 1),$$
where $(T_i)_{i \in \mathbb{N}}$ is defined by \eqref{iidme} and $N$ is the stopping time as in \eqref{Nstop}. Thus we have achieved the distribution of the three-dimensional Bessel process in Brownian motion. The existence of Brownian co-meander is obtained in the same vein.  $\square$ 
\begin{remark}
{\em Note that the stopping time $N$ defined as in \eqref{Nstop} has infinite mean, since $$\mathbb{E}N= \sum_{i \in \mathbb{N}} \mu^m_i(1)=\infty.$$ The problem whether Brownian co-meander or the three-dimensional Bessel process can be embedded in finite expected time, remains open. More generally, Rost \cite{Rost2} was able to characterize all stopping distributions of a continuous-time Markov process, given its initial distribution. In our setting, let $(P_t)_{t \geq 0}$ be the semi-group of the moving window process $X_t:=(B_{t+u}-B_t; 0 \leq u \leq 1)$ for $t \geq 0$, and $\mu^W$ be its initial distribution, corresponding to Wiener measure on $\mathcal{C}[0,1]$. Following Rost \cite{Rost2}, for any distribution $\mu$ on $\mathcal{C}[0,1]$, one can construct the continuous-time filling measures $(\mu_t,\mu^W_t)_{t \geq 0}$ and a suitable stopping time $T$ such that $$\mu-\mu_t+\mu^W_t=\mu^W P_{t \wedge T}.$$
Thus, the distribution $\mu$ is achieved if and only if $\mu_{\infty}=0$, where $\mu_{\infty}:=\downarrow \lim_{t \rightarrow \infty} \mu_t $. In particular, Brownian motion with drift $(\vartheta t+B_t; 0 \leq t \leq 1)$ for a fixed $\vartheta$, can be obtained for a suitable stopping time $T+1$.}
\end{remark}
\section{Potential theory for continuous-time patterns}
\label{4}
\quad In Question \ref{Q2}, we ask for any Borel measurable subset $S$ of $\mathcal{C}[0,1]$ whether $S$ is hit by the moving-window process $X_t:=(B_{t+u}-B_t;  0 \leq u \leq 1$ for $t \geq 0$, at some random time $T$. Related studies of the moving window process appear in several contexts. Knight \cite{Knight1, Knight2} introduced the prediction processes, where the whole past of the underlying process is tracked to anticipate its future behavior. The relation between  Knight's prediction processes and our problems is discussed briefly at the end of the section. Similar ideas are found in stochastic control theory, where certain path-dependent stochastic differential equations were investigated, see e.g. the monograph of Mohammed \cite{Mohammed} and Chang et al \cite{Chang}. More recently, Dupire \cite{Dupire} worked out a functional version of It\^{o}'s calculus, in which the underlying process is path-valued and notions as time-derivative and space-derivative with respect to a path, are proposed. We refer readers to the thesis of Fourni\'{e} \cite{Fournie} as well as Cont and Fourni\'{e} \cite{CF1,CF2,CF3} for further development.

\quad Indeed, Question \ref{Q2} is some issue of potential theory. In Benjamini et al \cite{BPP} a potential theory was developed for transient Markov chains on any countable state space $E$. They showed that the probability for a transient chain to ever visit a given subset $S \subset E$, is estimated by $\capacity_M(S)$ -- the Martin capacity of the set $S$. See also M\"{o}rters and Peres \cite[Section $8.3$]{PeresMorters} for a detailed exposition. As pointed out by Steven Evans \cite{StevePer}, such a framework still works well for our discrete renewal patterns. For $0 < \alpha<1$, define the $\alpha-$potential of the discrete patterns/strings of length $n$ as
\begin{align*}
G^{\alpha}(\epsilon',\epsilon'') :&=\sum_{k=0}^{\infty} \alpha^k P^{k}(\epsilon',\epsilon'') \\
&=\sum_{k=0}^{n-1}\left(\frac{\alpha}{2}\right)^k1\{\sigma_k(\epsilon')=\tau_k(\epsilon'')\}+ \frac{1}{1-\alpha} \left(\frac{\alpha}{2}\right)^k,
\end{align*}
where $\epsilon', \epsilon'' \in \{-1,1\}^n$, and $P(\cdot,\cdot)$ is the transition kernel of discrete patterns/strings of length $n$ in simple random walks, and $\sigma_k$ (resp. $\tau_k$): $\{-1,1\}^n \rightarrow \{-1,1\}^{n-k}$ the restriction operator to the last $n-k$ strings (resp. to the first $n-k$ strings). The following result is a direct consequence of the first/second moment method, and we leave the detail to readers.
\begin{proposition}\cite{StevePer}
Let $T$ be an independent $\mathbb{N}-$valued random variable with $\mathbb{P}(T>n)=\alpha^n$. For $\mathcal{A}^n$ a collection of discrete patterns of length $n$, we have
$$\frac{1}{2}\frac{2^{-n}}{1-\alpha}\capacity_{\alpha}(\mathcal{A}^n) \leq \mathbb{P}(T(\mathcal{A}^n)<T) \leq \frac{2^{-n}}{1-\alpha}\capacity_{\alpha}(\mathcal{A}^n) ,$$
where for $A \subset \{-1,1\}^n$,
$$\capacity_{\alpha}(A):=\left[\inf\left\{\sum_{\epsilon',\epsilon'' \in \{-1,1\}^n} G^{\alpha}(\epsilon',\epsilon'')g(\epsilon')g(\epsilon''); g \geq 0, g(A^c)=\{0\}~\mbox{and}\sum_{\epsilon \in \{0,1\}^n}g(\epsilon)=1\right\}\right]^{-1}.$$
\end{proposition}

\quad Now let us mention some previous work regarding the potential theory for path-valued Markov processes. There has been much interest in developing a potential theory for the Ornstein-Uhlenbeck process in Wiener space $\mathcal{C}[0,\infty)$, defined as
$$Z_t:=U(t, \cdot) \quad \mbox{for}~t \geq 0,$$
where $U(t,\cdot):=e^{-t/2}W(e^t,\cdot)$ is the Ornstein-Uhlenbeck Brownian sheet. Note that the continuous-time process $(Z_t; t \geq 0)$ takes values in the Wiener space $\mathcal{C}[0,\infty)$ and starts at $Z_0:=W(1,\cdot)$ as standard Brownian motion. Following Williams \cite{Williamsapp}, a Borel measurable set $S \subset \mathcal{C}[0,\infty)$ is said to be quasi-sure if $\mathbb{P}(\forall t \geq 0, Z_t \in S)=1$, which is known to be equivalent to
\begin{equation}
\label{qs}
\capacity_{OU}(S^c) =0, 
\end{equation}
where 
\begin{equation}
\label{qscap}
\capacity_{OU}(S^c):=\int_0^{\infty}e^{-t} \mathbb{P}(\exists T \in [0,t]~\mbox{such that}~Z_T \in S^c)dt
\end{equation}
is the {\em Fukushima-Malliavin capacity} of $S^c$, that is the probability that $Z$ hits $S^c$ before an independent exponential random time with parameter $1$. Taking advantage of the well-known {\em Wiener-It\^{o} decomposition} of the Ornstein-Uhlenbeck semigroup, Fukushima \cite{Fuku} provided an alternative construction of \eqref{qs} via the Dirichlet form. The approach allows the strengthening of many Brownian almost sure properties to quasi-sure properties. See also the survey of Khoshnevisan \cite{Khos} for recent development. 

\quad Note that the definition \eqref{qscap} can be extended to any (path-valued) Markov process. Within this framework, a related problem to Question \ref{Q2} is

{\bf Question ${\bf 1.5'}$:} {\em Given a Borel measurable set $S_{\infty} \subset \mathcal{C}[0,\infty)$, is 
\begin{align*}
\capacity_{MW}(S_{\infty}):&=\int_0^{\infty}e^{-t} \mathbb{P}[\exists T \in [0,t]~\mbox{such that}~ \Theta_T \circ B \in S_{\infty}]dt \\
&=0~\mbox{or }>0?
\end{align*}
where $(\Theta_t)_{t\geq 0}$ is the family of spacetime shift operators defined as
\begin{equation}
\label{shift}
\Theta_t \circ B:= (B_{t+u}-B_t; u \geq 0) \quad \mbox{for all}~t \geq 0.
\end{equation}}

\quad It is not difficult to see that the set function $\capacity_{MW}$ is a Choquet capacity associated to the shifted process $(B_{t+u}-B_t; u \geq 0)$ for $t \geq 0$, or the moving-window process $X_t:=(B_{t+u}-B_t; 0 \leq u \leq 1)$ for $t \geq 0$. For $S$ a Borel measurable subset of $\mathcal{C}[0,1]$, if $\capacity_{MW}(S \otimes_1 \mathcal{C}[0,\infty))=0$, where 
\begin{equation}
\label{concat}
S \otimes_1 \mathcal{C}[0,\infty):=\{(w_{t}1_{t< 1}+(w_1+w'_{t})1_{t \geq 1})_{t \geq 0}; w \in S~\mbox{and}~w' \in \mathcal{C}[0,\infty)\}
\end{equation}
is the usual path-concatenation, then $$\mathbb{P}(\exists T>0~\mbox{such that}~X_T \in S]=0,$$ i.e. almost surely the set $S$ is not hit by the moving-window process $X$. Otherwise, $$\mathbb{P}[\exists T \in [0,t]~\mbox{such that}~X_T \in S]>0 \quad \mbox{for some}~t \geq 0,$$
and an elementary argument leads to $\mathbb{P}[\exists T \geq 0~\mbox{such that}~X_T \in S]=1$.

\quad As context for this question, we note that path-valued Markov processes have also been extensively investigated in the superprocess literature. In particular, Le Gall \cite{LeGallp} characterized the polar sets for the Brownian snake, which relies on earlier work on the potential theory of symmetric Markov processes by Fitzsimmons and Getoor \cite{FitG} among others.

\quad There has been much progress in the development of potential theory for symmetric path-valued Markov processes. However, the shifted process, or the moving-window process, is not time-reversible and the transition kernel is more complicated than that of the Ornstein-Uhlenbeck process in Wiener space. So working with a non-symmetric Dirichlet form, see e.g. the monograph of Ma and R\"{o}ckner \cite{MaRo}, seems to be far from obvious.

\begin{openpb}~
\begin{enumerate}
\item
Is there any relation between the two capacities $\capacity_{X}$ and $\capacity_{MW}$ on Wiener space?
\item
Propose a non-symmetric Dirichlet form for the shifted process $(\Theta_t \circ B)_{t \geq 0}$, which permits to compute the capacities of the sets of paths $\mathcal{E}$, $\mathcal{M}$, $\mathcal{BR}^{\lambda}\ldots$etc. 
\end{enumerate}
\end{openpb}

\quad This problem seems substantial already for one-dimensional Brownian motion. But it could of course be posed also for higher dimensional Brownian motion, or a still more general Markov process. Following are some well-known examples of non-existing patterns in $d$-dimensional Brownian motion for $d \ge 2$.
\begin{itemize}
\item
$d=2$ (Evans \cite{Evans}): There is no random time $T$ such that $(B_{T+u}-B_T; 0 \leq u \leq 1)$ has a two-sided cone point with angle $\alpha<\pi$;
\item
$d=3$ (Dvoretzky et al \cite{DEKT}): There is no random time $T$ such that $(B_{T+u}-B_T; 0 \leq u \leq 1)$ contains a triple point;
\item
$d \geq 4$ (Kakutani \cite{Kaku}, Dvoretzky et al \cite{DEK}): There is no random time $T$ such that $(B_{T+u}-B_T; 0 \leq u \leq 1)$ contains a double point.
\end{itemize}
We refer readers to the book of M\"{o}rters and Peres \cite[Chapter $9$ and $10$]{PeresMorters} for historical notes and further discussions on sample path properties of Brownian motion in all dimensions.

\quad Finally, we make some connections between Knight's prediction processes and our problems. For background, readers are invited to Knight \cite{Knight1,Knight2} as well as the commentary of Meyer \cite{Meyer} on Knight's work. To avoid heavy measure theoretic discussion, we restrict ourselves to the classical Wiener space $(\mathcal{C}[0,\infty), \mathcal{F}, (\mathcal{F}_t)_{t \geq 0}, \mathbb{P}^{\bf W})$, where $(\mathcal{F}_t)_{t \geq 0}$ is the augmented Brownian filtrations satisfying the usual hypothesis of right-continuity. 

\quad The prediction process is defined as, for all $t \geq 0$ and $S_{\infty}$ a Borel measurable set of $\mathcal{C}[0,\infty)$,
$$Z^{{\bf W}}_t(S_{\infty}):=\mathbb{P}^{\bf W}[ \Theta_t \circ B \in S_{\infty} | \mathcal{F}_t],$$
where $\Theta_t \circ B$ is the shifted path defined as in \eqref{shift}. Note that $(Z^{{\bf W}}_t)_{t \geq 0}$ is a strong Markov process, which takes values in the space of probability measure on the Wiener space $(\mathcal{C}[0,\infty),\mathcal{F})$. In terms of the prediction process, Question \ref{Q2} can be reformulated as

{\bf Question ${\bf 1.5''}$:} {\em Given a Borel measurable set $S \subset \mathcal{C}[0,1]$, can we find a random time $T$ such that
$$\mathbb{E}Z_T^{{\bf W}}(S \otimes_1 \mathcal{C}[0,\infty))=1?$$
where $S \otimes_1 \mathcal{C}[0,\infty)$ is defined as in \eqref{concat}.}\\\\
\textbf{Acknowledgement:} We would like to express our gratitude to Patrick Fitzsimmons for posing the question whether one can find the distribution of Vervaat bridges by a random spacetime shift of Brownian motion. We thank Steven Evans for helpful discussion on potential theory, and Davar Koshnevisan for remarks on additive L\'{e}vy processes.
\bibliographystyle{plain}
\bibliography{Patterns}
\end{document}